\documentclass[10pt]{amsart}
\usepackage[margin=1.5in]{geometry}
\usepackage{amssymb,amsthm,amsmath}
\usepackage[utf8]{inputenc}
\usepackage{preamble}

\title{$K$-tropical curves, Riemann-Roch and Serre Duality}
\author{Soham Karwa, Joe Rabinoff}
\begin{document}
\begin{abstract}
    We endow a tropical curve with the structure of a ringed space and establish analogues of Riemann-Roch and Serre duality for these ringed spaces. As a corollary, we give a cohomological proof of the tropical Riemann-Roch theorem. 
\end{abstract}
\maketitle
\vspace{-2em}

\section{Introduction}

\subsection{Main results}
The Riemann--Roch theorem for algebraic curves relates the degree of a divisor $D$ to the dimension of the space of meromorphic functions with zeros and poles constrained by $D$. Remarkably, an analogous formula holds for tropical curves, with the dimension of a linear system replaced by tropical divisor rank. Baker and Norine \cite{bakernorine} first established the corresponding theorem for finite graphs using chip-firing and reduced divisors; the extension to metric graphs was then obtained independently by Gathmann and Kerber \cite{gathmann} and Mikhalkin and Zharkov \cite{mikhalkinzharkov}.

A longstanding goal in tropical geometry has been to give a cohomological proof of the tropical Riemann--Roch theorem. In this paper, we introduce $K$-tropical curves, obtained by equipping tropical curves (\ie metric graphs) with sheaves of $K$-algebras arising from $K$-affine structures \cite{ksaffine, karwa2025nonarchimedeanperiodslogcalabiyau}, where $K$ is a complete non-archimedean field of residue characteristic zero. We establish Serre duality and Riemann--Roch for these ringed spaces and relate their cohomology to tropical divisor rank, thereby obtaining such a proof. We will denote a tropical curve by $\underline{\sqsubset}$ in this paper.
\begin{thm}[$K$-tropical Serre Duality]
            \label{ktropserre}
    Let $p_\sqsubset: (\underline{\sqsubset},\Oc_\sqsubset^\an) \rightarrow \spec K$ be a $K$-tropical curve and $\Lc$ a line bundle on $(\underline{\sqsubset},\Oc_\sqsubset^\an)$. Then there exists a canonical line bundle $\omega_{\sqsubset^\an}$ on $(\underline{\sqsubset},\Oc_{\sqsubset}^\an)$ such that the cohomology groups of $\Lc$ are finite dimensional and  
    \begin{align}
        H^i(\Lc)^\vee \cong H^{1-i}( \omega_{\sqsubset^\an}\otimes \Lc^{-1}).
    \end{align}
\end{thm}
As is done in the classical case, we use $K$-tropical Serre duality to prove the $K$-tropical Riemann-Roch theorem by first proving a weaker version. \begin{thm}
    \label{cheaprr}
    Let $\Lc$ be a line bundle on $(\underline{\sqsubset},\Oc_\sqsubset^\an)$. Then \begin{align}
        \chi(\Lc) = \chi(\Oc_\sqsubset^\an)  + \deg \Lc.
    \end{align}
\end{thm}

Suppose the genus of $\underline{\sqsubset}$ is $g$. Then an explicit computation shows $\chi(\Oc_\sqsubset^\an) = 1-g$ which together with $K$-tropical Serre duality gives the $K$-tropical Riemann-Roch theorem. 

\begin{thm}[$K$-tropical Riemann-Roch]
    \label{thm:Ktroprr}
    Let $\Lc$ be a line bundle on $(\underline{\sqsubset},\Oc_\sqsubset^\an)$. Then \begin{equation}
        h^0(\Lc) - h^0(\omega_{\sqsubset^\an}\otimes \Lc^{-1}) = \deg \Lc -g +1. 
    \end{equation}
\end{thm}
We establish a divisor-line bundle correspondence for $K$-tropical curves and give an analogue of Baker's specialization lemma \cite{specialbaker} for $K$-tropical curves relating the $K$-tropical rank, defined by $h^0(\Lc)-1$, to the tropical rank of a divisor. To deduce the tropical Riemann-Roch theorem, we take $K = \C((t^\R))$ to be the field of Puiseaux series over $\C$ and enhance $\underline{\sqsubset}$ with the structure of a $K$-tropical curve induced from the lift of the tropical period matrix to a period matrix over $K$. In this case, we prove the reverse inequality of the specialization lemma and thus establish equality of the $K$-tropical rank and tropical rank of a divisor. 
\begin{cor}[Tropical Riemann-Roch]
Let $\underline{\sqsubset}$ be a genus $g$ tropical curve and $\omega_\sqsubset^\trop$ the canonical tropical line bundle on $\underline{\sqsubset}$. Let $\Lc$ be a a tropical line bundle on $\underline{\sqsubset}$. Then
\begin{equation}
    r_{\sqsubset,\trop}(\Lc) - r_{\sqsubset,\trop}(\omega_\sqsubset^\trop\otimes \Lc^{-1}) = \deg \Lc - g +1. 
\end{equation}
where $r_{\sqsubset,\trop}(\cdot)$ is the tropical divisor rank. 
\end{cor}

Given a finite connected graph $G$ without loop edges, let $\underline{\sqsubset}$ be the corresponding metric graph in which every edge has length one. Then the divisor rank on the graph $G$ and the metric graph $\underline{\sqsubset}$ are equal \cite{luo,hladkykravnorine} and thus we can deduce the Riemann-Roch theorem for graphs. 

\begin{cor}[Riemann-Roch for graphs]
    Let $G$ be a finite connected genus $g$ graph without loop edges. Let $r_G$ denote the divisor rank on $G$ and $K_G$ the canonical divisor on $G$. Given a divisor $D$ on $G$, then \begin{align*}
        r_G(D) - r_G(K_G-D)= \deg D - g + 1.
    \end{align*}
\end{cor}

The $K$-tropical Riemann-Roch theorem can also be used to deduce the algebraic Riemann-Roch theorem for maximally degenerate curves. Recall, a proper curve $C$ over $K$ is maximally degenerate if, after a finite extension of $K$, it has totally split degenerate reduction, \ie the special fiber $\mathscr{C}_k$ of its minimal snc model $\mathscr{C}$ is reduced and all its irreducible components are rational. Then the essential skeleton of $C$, which can be described as the dual intersection complex of $\mathscr{C}_k$, embeds into the Berkovich analytification $C^\an$ of $C$. Moreover, it has the structure of a tropical curve $\underline{\sqsubset}$ and there is a retraction map $\rho: C^\an \longrightarrow \underline{\sqsubset}$. 
 Equipping $\underline{\sqsubset}$ with the sheaf of rings $\rho_*(\Oc_{C^\an})$, we can deduce the algebraic Riemann-Roch theorem from the $K$-tropical Riemann-Roch theorem.

\begin{cor}[Algebraic Riemann-Roch]
    Let $C$ be a maximally degenerate genus $g$ curve and $\Lc$ a line bundle on $C$. Then \begin{align}
        h^0(\Lc) - h^0(\omega_C\otimes \Lc^{-1}) = \deg \Lc - g + 1. 
    \end{align}
\end{cor}

Whilst every tropical curve can be realized as the essential skeleton of a maximally degenerate tropical curve, it is well known that the algebraic Riemann-Roch theorem cannot na\"ively imply the tropical Riemann-Roch theorem and vice versa as there are examples where given a tropical curve and a divisor, there is no algebraic lift that realizes the tropical rank \cite[Example 5.13]{liftharmonic}. The ringed space structure we construct to prove tropical Riemann-Roch does not necessarily arise from an algebraic curve. Indeed, lifting the tropical period matrix to a period matrix over $K$ corresponds to constructing a maximally degenerate abelian variety over $K$ whose essential skeleton is the tropical  Jacobian of $\underline{\sqsubset}$. This abelian variety may not lie in the Torelli locus and thus need not arise from the Jacobian of an algebraic curve. 

\begin{rmk}
    In deducing the algebraic Riemann-Roch theorem, the ringed space structure on $\underline{\sqsubset}$ does not arise from a $K$-affine structure in the current paper. We nevertheless refer to it as a $K$-tropical curve but will distinguish between the two cases. In fact, in Example \ref{ex:noKaffinealg}, we prove that such a $K$-affine structure cannot exist if we also require that the Picard group of $C$ is isomorphic to the Picard group of the associated $K$-tropical curve $\sqsubset$. In future work, we will extend the notion of $K$-affine structure to give a more universal construction of a $K$-tropical curve. 
\end{rmk}

\subsection{Strategy}

We outline the construction of the ringed space structure on $\underline{\sqsubset}$ used in the proof of the tropical Riemann-Roch theorem and the proofs of several of the results above. Let $\Aff_{\sqsubset,\Z}$ denote the sheaf of integral affine functions on $\underline{\sqsubset}$. This fits in the short exact sequence 
\begin{align*}
    0 \longrightarrow \R \longrightarrow \Aff_{\Z,\sqsubset} \longrightarrow \check{\Lambda}_\sqsubset \longrightarrow 0,
\end{align*}
where $\Aff_{\Z,\sqsubset}$ is the sheaf of integral affine functions on $\underline{\sqsubset}$. In the proof of the tropical Riemann-Roch theorem, we will take $K= \C((t^\R))$ so that we can also deal with tropical curves with real edge lengths. By picking the section $s: \R \longrightarrow \C((t^\R)), a \mapsto t^a$ and using covariance of $\text{Ext}^1$ in the second entry, we have a map $\text{Ext}^1(\check{\Lambda}_\sqsubset,\R) \longrightarrow \text{Ext}^1( \check{\Lambda}_\sqsubset, K^\times)$. We take the $K$-affine structure $\Aff_{K,\sqsubset}$ to be the image of $\Aff_{\Z,\sqsubset}$. We describe $\Aff_{K,\sqsubset}$ explicitly in \S\ref{kaffinestruc}. By construction, we then have a diagram 
\begin{equation}
\begin{tikzcd}[column sep=large, row sep=large]
0 \arrow[r] &
K^\times
\arrow[r]
\arrow[d, "\operatorname{val}"]
&
\operatorname{Aff}_{K,\sqsubset}
\arrow[r]
\arrow[d]
&
\check{\Lambda}_{\sqsubset}
\arrow[r]
\arrow[d, "\mathrm{id}"]
&
0
\\
0 \arrow[r] &
\mathbb{R}
\arrow[r]
\arrow[u, bend left=35]
&
\operatorname{Aff}_{\mathbb{Z},\sqsubset}
\arrow[r]
\arrow[u, bend left=35]
&
\check{\Lambda}_{\sqsubset}\arrow[u, bend left=35]
\arrow[r]
&
0
\end{tikzcd}
\end{equation}

    In \S\ref{ktropsheaf}, we define the ringed space structure to be a suitable quotient of the group algebra $K[\Aff_{K,\sqsubset}]$. We denote the sheaf of rings by $\Oc_\sqsubset$ and the ringed space by $\sqsubset = (\underline{\sqsubset},\Oc_\sqsubset)$. We can also define the ringed space structure on the tropical Jacobian analogously via the $K$-affine structure and the section $s$.  By functoriality, the tropical Abel map is then upgraded to a morphism of ringed spaces.
 	
 	In \S\ref{analyticmod}, by an appropriate completion, we  modify the structure sheaf $\Oc_\sqsubset$ to produce $\Oc_\sqsubset^\an$ which was used in the statements of the theorems above. This is necessary for several reasons including a sensible notion of Cartier divisor and more importantly the existence of global sections. We then turn our attention to study the line bundles on the $K$-tropical curve $(\underline{\sqsubset}, \Oc_\sqsubset^\an)$ where we explicitly describe the Picard group and define the degree of a line bundle. We define Weil divisors on $\sqsubset$ by identifying the appropriate subset of $K$-points which arise from line bundles. Using the $d$-fold $K$-tropical Abel map, we define linear equivalence of divisors.  This is carried out in \S \ref{ktropsheaf} and \S\ref{ktropdivisor}. In \S \ref{divisorlinebundle}, Cartier divisors are introduced and a divisor-line bundle correspondence is proved. 
 	
 	In \S\ref{globalsections}, we define the linear system of a $K$-tropical divisor $D$ and its rank. We show the rank of an effective divisor $D$ is equal to $h^0(\Oc_{\sqsubset^\an}(D))-1$ where $\Oc_{\sqsubset^\an}(D)$ is the corresponding line bundle. In particular, we show that the cohomology of a line bundle on our ringed space is finite dimensional. In this section, we also give a $K$-tropical analogue of the specialization lemma of Baker and end the section by proving Theorem \ref{cheaprr}.

 In \S \ref{serreduality}, we prove Serre duality by appealing to a more general result of \cite{spaltenstein} concerning a duality theorem for ringed spaces. Indeed, we show that the the dualizing sheaf $p_\sqsubset^! K$ is a line bundle concentrated in degree one, which we denoted $\omega_{\sqsubset^\an}$ in Theorem \ref{ktropserre}. Using the finite dimensionality of cohomology of line bundles, we then prove Theorem \ref{ktropserre}. With Theorem \ref{cheaprr} and an explicit computation that $\chi(\underline{\sqsubset},\Oc_\sqsubset^\an) = 1-g$, we prove the $K$-tropical Riemann-Roch theorem. As a corollary, we deduce the algebraic Riemann-Roch theorem.
 
 In the final section  \S\ref{troprr},  given a tropical curve $\underline{\sqsubset}$, we explain how the constructions in the previous sections give a cohomological proof of the tropical Riemann-Roch theorem. 
  	
\subsection{Future directions}

We have primarily concerned ourselves with line bundles in this paper. A natural next step would be to consider vector bundles and coherent sheaves. One reason for this is our methods do not immediately extend to this case as we use the $d$-fold $K$-tropical Abel map to prove finite dimensionality of cohomology groups. Nevertheless, Serre duality currently does give a duality between the relevant cohomology groups. Once finite dimensionality is established, it would be interesting to prove a Grothendieck-Riemann-Roch formula for vector bundles. Furthermore, since the techniques of this paper allow one to use cohomological methods in tropical geometry, a potential direction to explore is the study of tropical Brill-Noether loci. 

An active area of interest is the correct notion of tropical vector bundles and the compactification of the moduli space of vector bundles on special fibers of maximally degenerate curves. The ringed space structure on the $K$-tropical curve may be useful here. These ideas will be explored in future work. 

\subsection{Acknowledgments}
The authors would like to thank Ravi Vakil for suggesting that the results of \cite{karwa2025nonarchimedeanperiodslogcalabiyau} could be relevant to proving the tropical Riemann-Roch theorem cohomologically. They would also like to thank Omid Amini for many helpful suggestions to realize the ultimate goal of proving tropical Riemann-Roch. 

\subsection{AI Disclosure}
The authors used ChatGPT-6 Astra for references and proof reading. The ideas and proofs are entirely the authors’ own.
\tableofcontents

\section{Preliminaries}


\subsection{Tropical curves and their Jacobians}

Throughout this paper, a tropical curve $\underline{\sqsubset}$ will be a finite connected metric graph with at least one edge. When it is clear from context, we will simply denote the tropical curve by $\sqsubset$. We will denote the vertex (resp. edge) set by $V(\sqsubset)$ (resp. $E(\sqsubset)$). Let the \textit{genus} of $\sqsubset$ be $g(\sqsubset) = b_1(\sqsubset)$. We will assume that all graphs are connected and there are no 1-valent vertices.   A \textit{bridge} of $\sqsubset$ is an edge of $\sqsubset$ which is not contained in any non-trivial cycle.  Given a vertex $v \in V(\sqsubset)$, let $\Star(v)\subseteq \sqsubset$  denote the open star of $v$. Two tropical curves $\sqsubset, \sqsubset'$ are equivalent if there is a common length preserving refinement. In particular, we will freely subdivide edges by inserting $2$-valent vertices whenever convenient. Such subdivisions do not alter the underlying metric graph or the associated tropical curve, but merely change the chosen combinatorial model. 

\begin{defn}(Harmonic function)
	Let $\Gamma$ be a finite metric graph and let $F: \Gamma \rightarrow \R$ be continuous function. Then 
	\begin{enumerate}
		\item We say that $F$ is piecewise affine with integer slopes if the restriction of $F$ to each edge (identified with an interval) has the form $F(x) = mx + b$ with $m \in \Z$ and $b \in \R$ and if there are only finitely many points where $F$ is not differentiable; 
		\item For $x \in \Gamma$, let $T_x$ denote the set of tangent directions at $x$ and for $v \in T_x$, we let $d_vF$ denote the outgoing slope of $F$ in the $v$-direction. We say that $F$ is \textit{harmonic} at $x \in \Gamma$ provided that $\sum_{v \in T_x} d_vF = 0$ and we say that $F$ is harmonic if it is harmonic at all points $x \in \Gamma$. 
	\end{enumerate}
	A tropical curve $\sqsubset$ determines a sheaf of harmonic functions which we denote $\Aff_{\sqsubset,\Z}$.  Note that subdivision of edges in $\sqsubset$ does not change $\Aff_{\sqsubset,\Z}$. 
\end{defn}

Following \cite[Definition 4.1]{mikhalkinzharkov}, the integral cotangent bundle $\check{\Lambda}_{\sqsubset}$ is defined via the short exact sequence  
	\begin{align}
		0 \rightarrow \R \rightarrow \Aff_{\sqsubset,\Z}\rightarrow \check{\Lambda}_\sqsubset \rightarrow 0
		\label{AffZ}
	\end{align}
where $\R$ is the sheaf of constant functions on $\sqsubset$. At a point $p \in \sqsubset$, the stalk $\check{\Lambda}_{\sqsubset,p}$ is given by $\{f \in \Aff_{\sqsubset,p}: f(p) = 0\}$. 

\subsubsection{Divisors on tropical curves}
A \textit{divisor} on $\sqsubset$ is a formal finite integer linear combination
$D=\sum_{p\in\sqsubset}a_p\cdot p , \ a_p\in\Z$,
where all but finitely many coefficients $a_p $vanish. The group of divisors is denoted $\operatorname{Div}(\sqsubset)$, and the degree of $D$ is
$\deg(D)=\sum_{p}a_p.$
A divisor is \textit{effective} if $a_p\ge0$ for every $p$. A divisor $D$ is \textit{principal} if there exists a continuous piecewise integral affine function $f:\ \sqsubset \longrightarrow \R$  such that \[D = \text{div}(f) = \sum_{p \in  \sqsubset}\text{ord}_p(f)\cdot p 
\]
where $\text{ord}_p(f)$ is the sum of all outgoing slopes at $p$.  
Two divisors $D,D'$ on $\sqsubset$  are equivalent if and only if $D-D' = \text{div}(f)$ for some continuous piecewise integral affine function $f:\ \sqsubset \longrightarrow \R$. We denote the tropical divisors on $\sqsubset$ by $\text{TroDiv}(\sqsubset)$ and denote the tropical divisors of degree $d$ by $\text{TroDiv}^d(\sqsubset)$. We denote the tropical principal divisors by $\text{TroPrin}(\sqsubset)$.

The following lemma allows us to write a $k$-valent vertex of $\sqsubset$ for $k\geq 3$, viewed as a degree 1 divisor, as a sum of neigbouring $2$-valent vertices. This will be a very useful operation when we introduce the ringed space structure on $\sqsubset$ as a means to describe effective Cartier divisors on $\sqsubset$ in \S\ref{globalsections}.

\begin{lem}
	Let $\sqsubset$ be a tropical curve and $v \in V(\sqsubset)$ be a $k\geq 3$ valent vertex. Then after adding $2$-valent vertices $q_i$ in $\Star(v)$, there is a tropical rational function $f$ such that \begin{align*}
		\text{div}(f) = v + \sum m_i q_i, m_i \in \Z
	\end{align*}
\end{lem}

\begin{proof}
Consider the star of $v$ as in Figure \ref{starv} and choose a branch $e_1$. There exists a rational number $1/q \in \Q$ such that a closed interval $I_{1/q}$ of length $1/q$ embeds into each branch of $\Star(v)$. Choose an embedding of $I_{1/q}$ for each branch such that on $e_1$, one of the endpoints of the intervals is $v$ and for all other branches the embedding is strict inside the branch. One can then define the piecewise linear function $f$ on $\Star(v)$ as follows: \begin{enumerate}
	\item at $v$, the value is 0; 
	\item on $e_1$, $f$ has slope 1 on the embeded $I_{1/q}$ and 0 elsewhere;
	\item on $e_i, \ i \neq 1$, $f$ has slope 1 on the embeded $I_{1/q}$ and 0 elsewhere;
	\item the function $f$ is constant outside $\Star(v)$. 
 \end{enumerate}
At the endpoints in the $\Star(v)$, $f$ takes the same value and thus the constant extension in (4) is possible. By construction, the divisor of $f$ is of the required form.

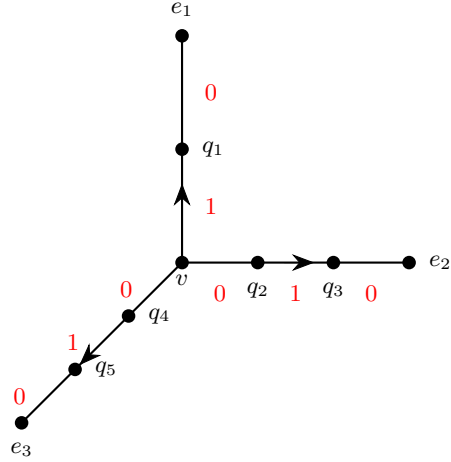
\begin{figure}[h!]
	
\begin{tikzpicture}[scale =0.5,
    >=Stealth,
    vertex/.style={circle,fill=black,inner sep=1.8pt},
    edge/.style={line width=0.8pt},
    edgelabel/.style={red,font=\small},
    every node/.style={font=\small}
]


\coordinate (v) at (0,0);

\coordinate (e1) at (0,6);
\coordinate (e2) at (6,0);
\coordinate (e3) at (-4.24,-4.24); 

\coordinate (q1) at ($(v)!0.5!(e1)$);

\coordinate (q2) at ($(v)!1/3!(e2)$);
\coordinate (q3) at ($(v)!2/3!(e2)$);

\coordinate (q4) at ($(v)!1/3!(e3)$);
\coordinate (q5) at ($(v)!2/3!(e3)$);


\draw[edge] (v)--(e1);
\draw[edge] (v)--(e2);
\draw[edge] (v)--(e3);


\foreach \x in {v,q1,q2,q3,q4,q5,e1,e2,e3}
    \node[vertex] (\x-node) at (\x) {};


\node[below =1pt] at (v) {$v$};

\node[right=4pt] at (q1) {$q_1$};
\node[below=4pt] at (q2) {$q_2$};
\node[below=4pt] at (q3) {$q_3$};
\node[right=4pt] at (q4) {$q_4$};
\node[right=4pt] at (q5) {$q_5$};

\node[above=4pt] at (e1) {$e_1$};
\node[right=4pt] at (e2) {$e_2$};
\node[below=4pt] at (e3) {$e_3$};


\draw[-{Stealth[length=3mm]}]
    ($(v)!0.30!(q1)$)--($(v)!0.70!(q1)$);

\draw[-{Stealth[length=3mm]}]
    ($(q2)!0.35!(q3)$)--($(q2)!0.75!(q3)$);

\draw[-{Stealth[length=3mm]}]
    ($(q4)!0.55!(q5)$)--($(q4)!0.92!(q5)$);


\node[edgelabel,right=5pt] at ($(q1)!0.5!(e1)$) {$0$};
\node[edgelabel,right=5pt] at ($(v)!0.5!(q1)$) {$1$};

\node[edgelabel,below=5pt] at ($(v)!0.5!(q2)$) {$0$};
\node[edgelabel,below=5pt] at ($(q2)!0.5!(q3)$) {$1$};
\node[edgelabel,below=5pt] at ($(q3)!0.5!(e2)$) {$0$};

\node[edgelabel,left=5pt] at ($(v)!0.5!(q4)$) {$0$};
\node[edgelabel,left=5pt] at ($(q4)!0.5!(q5)$) {$1$};
\node[edgelabel,left=5pt] at ($(q5)!0.5!(e3)$) {$0$};

\end{tikzpicture}
\caption{The star of $v$, the slope of $f$ is depicted in red}
\label{starv}
\end{figure}

\end{proof}


\begin{cor}
	Any $k$-valent vertex is linearly equivalent to a sum of $2$-valent vertices. 
	\label{kvalent2valent}
\end{cor}

\subsubsection{Tropical line bundles}

Whilst tropical divisors are usually easier to work with in practice, there is a notion of tropical line bundle as an $\Aff_{\sqsubset,\Z}$-torsor on $\sqsubset$. When we introduce a ringed space structure on $\sqsubset$, line bundles will be easy to define whereas we will need further machinery to make sense of a divisor on these spaces. 
\begin{defn}
	The tropical Picard group $\text{TroPic}(\sqsubset)$ is defined as  $H^1(\sqsubset,\Aff_{\sqsubset,\Z})$ and classifies isomorphism classes of tropical line bundles on $\sqsubset$. 
\end{defn}
By \cite[Corollary 4.7]{mikhalkinzharkov}, there is a divisor-line bundle correspondence \[\text{TroPic}(\sqsubset) \cong \text{TroDiv}(\sqsubset)/\text{TroPrin}(\sqsubset).\] Given a tropical divisor $D$, let $\Oc_{\sqsubset,\trop}(D)$ denote the corresponding line bundle.

\subsubsection{Tropical Jacobians}


\begin{definition}
A \textit{principally polarized tropical abelian variety} is a triple $(A, \Lambda, Q)$, where $A$ is a $g$-dimensional real torus, $\Lambda \cong \mathbb Z^g$ is a free abelian group of rank $g$, and $Q\colon \Lambda\times \Lambda\to \R$ is a symmetric, positive-definite bilinear form, such that we are given an isomorphism $A \cong \mathrm{Hom}(\Lambda, \mathbb R) / Q(\Lambda)$, where we write $Q : \Lambda \to \mathrm{Hom}(\Lambda, \mathbb R)$ for the natural map derived from the bilinear form.
\end{definition}

The tangent space $T_pA$ at a point $p \in A = V/\Lambda$ can be identified with $V$ and $\Lambda'$ can be viewed as a lattice of sections in $TA. $ An element $ \lambda \in \Lambda'$ is called an \textit{integral tangent vector}. 

 Let $\Omega^1_\sqsubset := \check{\Lambda}_\sqsubset(\sqsubset)$ denote the space of global 1-forms on $\sqsubset$. 
 An element $\gamma \in H_1(\sqsubset,\Z)$ can be viewed as an element of $\Omega^{1,\vee}_\sqsubset{\otimes\mathbb R} = \Hom(\Omega_{\sqsubset}^1, {\mathbb R})$ via the integration pairing $
\mu: \gamma \mapsto \int_\gamma$. Then the \textit{tropical Jacobian} is defined as \begin{align*}
    \text{TroJac}(\sqsubset) := \Omega^{1,\vee}_\sqsubset\otimes \R/\mu(H_1(\sqsubset,\Z)).
\end{align*}
Setting $M = H_1(\sqsubset,\Z)$ and identifying $\Omega_\sqsubset^1$ with $H_1(\sqsubset,\Z)$ in the obvious way, then the above definition becomes 
\begin{align*}
    \text{TroJac}(\sqsubset) := N_\R /\mu(H_1(\sqsubset,\Z))
\end{align*}
where $\mu: H_1(\sqsubset,\Z) \rightarrow N_\R := \Hom(M,\R)$ is given by the edge length pairing \cite[(2.7.1)]{jacobian}.  The edge length pairing $\langle \cdot, \cdot \rangle : M \times M \rightarrow \R$ makes $\text{TroJac}(\sqsubset)$ a principally polarized tropical abelian variety whose underlying topological space is homeomorphic to a $g$-dimensional real torus. There is a sheaf of integral affine functions $\Aff_{\text{TroJac}(\sqsubset),\Z}$ on $\text{TroJac}(\sqsubset)$  which fits in the short exact sequence 
\begin{align}
0 \longrightarrow \R \longrightarrow\Aff_{\text{TroJac}(\sqsubset) } \longrightarrow \check{\Lambda}_{\text{TroJac}(\sqsubset)} \longrightarrow 0
\end{align}
 with $\check{\Lambda}_{\text{TroJac}(\sqsubset)} \cong M$. 

Note $H^0(\sqsubset, \R) \cong \R \cong H^0(\Aff_{\sqsubset,\Z})$ since $H^0(\sqsubset,\R) \hookrightarrow H^0(\sqsubset,\Aff_{\sqsubset,\Z})$ embeds constant functions into harmonic functions but as $\sqsubset$ is compact and connected, $H^0(\Aff_{\sqsubset,\Z}) \cong \R$. Hence taking long exact sequences of (\ref{AffZ}), gives \begin{align}
	0 \rightarrow \Omega^1_\sqsubset \rightarrow H^1(\sqsubset,\R) \rightarrow \text{TroPic}(\sqsubset) \rightarrow H^1(\sqsubset,\check{\Lambda}_\sqsubset) \rightarrow 0.
\end{align}
By \cite[9.3.16]{mikhalkin2009tropical}, $H^1(\sqsubset,\check{\Lambda}_\sqsubset) \cong \Z$ and the map $\text{TroPic}(\sqsubset) \rightarrow H^1(\sqsubset, \check{\Lambda}_\sqsubset) $ coincides with the degree map and $\Omega_\sqsubset^1 \rightarrow H^1(\sqsubset,\R)$ is the map given by integration. We denote the degree $d$ component of $\TroPic(\sqsubset)$ by $\TroPic^d(\sqsubset)$.  In particular, this implies the tropical Abel-Jacobi theorem \begin{align*}
	\text{TroPic}^0(\sqsubset) \cong \text{TroJac}(\sqsubset). 
\end{align*}

\subsubsection{Linear systems and rank}

Given a tropical divisor $D$ on $\sqsubset$, one can associate the \textit{tropical linear system} \begin{align*}
	|D| := \{ E \in \text{Div}(\sqsubset) : E\sim D, E \geq 0\}.
\end{align*}
The tropical linear system $|D|$ has the structure of a polyhederal complex which may not be equidimensional. In particular, care is required when  defining the rank. For example, if one defines the rank as the maximal dimensional cell of $|D|$, it can be shown that there is no version of the tropical Riemann-Roch theorem \cite[Remark 8.4.2]{mikhalkin2009tropical}.  To remedy this, Baker-Norine \cite{bakernorine} introduced the notion of combinatorial rank defined as follows \begin{align*}
	r_{\sqsubset,\trop}(D) = \max \{ n : |D-E| \neq \varnothing \text{ for all effective divisors } E\in \text{Div}^n(\sqsubset)\}.
\end{align*}
where $r_{\sqsubset,\trop}(D) = -1 $ if $|D| = \emptyset$. The combinatorial rank satisfies the Riemann-Roch formula \begin{align}
	r_{\sqsubset,\trop}(D) - r_{\sqsubset,\trop}(K_\sqsubset^\trop-D) = \deg D- g +1.
\end{align}
where $K_\sqsubset^\trop := \sum_{v \in V(\sqsubset)} (\val(v) -2) \cdot v$ is the \textit{tropical canonical divisor}. We record the following straightforward lemma.

\begin{lem}
    Let $D$ be a degree $0$ tropical divisor on $\sqsubset$.  Then \begin{align*}
        r_{\sqsubset,\trop}(D) = \begin{cases}
            0 \  \ &\text{if } D \sim 0 \\ -1 &\text{otherwise }
        \end{cases}
    \end{align*}
    \label{degree0trop}
\end{lem}









\subsubsection{The tropical Abel-Jacobi map}

\begin{defn}
    Given a tropical abelian variety $(A = V/\Lambda, \Lambda')$ and a tropical curve $\sqsubset$, a continuous map $ f: \sqsubset \rightarrow A$ is \textit{tropical} if the following three conditions hold 
    \begin{enumerate}
    \item $f$ is piecewise affine 
    \item at each point $p \in \sqsubset$, the differential map $df_p: T_p\sqsubset \longrightarrow T_{f(p)}A$ sends integral tangent vectors to integral tangent vectors 
    \item at each point $p \in \sqsubset$, $f$ is harmonic \ie $\sum_\lambda df_p(\lambda) = 0$ where the sum is over all primitive tangent vectors at $p$. 
    \end{enumerate}
\end{defn}

\begin{defn}
    Given a vertex $p \in \sqsubset$, the morphism 
\begin{align*}
    \alpha_p: \sqsubset &\longrightarrow \TroJac(\sqsubset), \\ 
    q &\longmapsto \Oc(q-p)
\end{align*}
is the called the \textit{tropical Abel-Jacobi map}. By \cite[Lemma 9.3.8]{mikhalkin2009tropical}, $\alpha_p$ is a tropical morphism and thus induces a map of sheaves \[\alpha_p^{-1}\check{\Lambda}_{\Jac(\sqsubset)} \longrightarrow \check{\Lambda}_{\sqsubset}.\]
One can also define the \emph{Abel map} $\alpha: \sqsubset\rightarrow \text{TroPic}^1(\sqsubset), q \mapsto \Oc(q)$ which is also a tropical morphism. 

\end{defn}

Just as in the case of algebraic curves, for every integer $d\geq 0$, there is the $d$-fold Abel map \begin{align}
 \alpha_d: \sqsubset_d \ \longrightarrow	\text{TroPic}^d(\sqsubset)
 \end{align}
 where $\sqsubset_d$ is the $d$th symmetric product of $\sqsubset$. We will denote the tropical Abel map $\alpha_1$ simply by $\alpha$. 
\begin{thm}
	The following are true:
	\begin{enumerate}
	\item The $d$-fold Abel map is also a tropical morphism 

		\item The map $\alpha_1$ contracts bridge edges and is an isometry elsewhere;
		\item The map $\alpha_d$ is surjective if and only if $d \geq g$. 
		\item If $D$ is a divisor of degree $d$, then  $\alpha_d^{-1}(\Oc(D)) =|D|$.
	\end{enumerate}
\end{thm}

\begin{proof}
	(1) follows by \cite[Proposition 4.1]{brandtulirsch}.  (2) follows by Theorem 1.8 and Lemma 4.6 in \cite{bakernorine} and noting $\alpha_d$ translation of the $d$-fold tropical Abel-Jacobi map (3) is Theorem 1.7 in \cite{bakernorine} and (4) follows by definition of the map $\alpha_d$. 

\end{proof}

\subsection{Tropical modification}
\label{tropmod}
Any graph can be decomposed into a union of its maximal 2-connected subgraphs and of its bridges. We shall write such a decomposition as \[\sqsubset = \left(\cup_i \sqsubset_i \right)\cup B\]where $B$ is the set of bridge edges. By \cite[Corollary 2.14]{bakerfaber}, the tropical Jacobian of $\sqsubset$ decomposes according to the maximal 2-connected subgraphs \ie $\TroJac(\sqsubset)$ is canonically isomorphic to $\TroJac(\sqsubset_1) \times \dots \times \TroJac(\sqsubset_k)$. A \textit{cut vertex} of a graph $G$ is a vertex $v$ such that $G\backslash \{v\}$ is disconnected; the maximal 2-connected subgraphs $\sqsubset_i$ do not have cut vertices. By construction of the tropical Abel map, the tropical Abel map $\alpha$ when restricted to $\sqsubset_i$ then factors through $\text{TroJac}(\sqsubset_i)$. 


 In our construction of the ringed space on a tropical curve, we will need to further analyze the geometry of the tropical Abel map. We begin with the following combinatorial lemma. 

\begin{lem}
Let $G$ be a connected graph of genus $g$, and let $v$ be a vertex of valence $d$. If $d> g+1$, then $G\setminus\{v\}$ is disconnected.
\label{disconnected}
\end{lem}

\begin{proof}
Let $c$ denote the number of connected components of $G\setminus\{v\}$. Since removing $v$ deletes $d$ edges and one vertex,
\[
|V(G\setminus\{v\})|=|V(G)|-1,\qquad
|E(G\setminus\{v\})|=|E(G)|-d.
\]
By the formula for the Betti genus of a graph, we have \begin{align*}
g(G\setminus\{v\})
&= |E(G\setminus\{v\})|
   -|V(G\setminus\{v\})|
   +c \\
&= g-d+c.
\end{align*}
Suppose, for contradiction, that $G\setminus\{v\}$ is connected. Then $c=1$, so
\[
g(G\setminus\{v\})=g-d+1.
\]
Since $d>g+1$, the right-hand side is negative, contradicting the fact that the genus of a graph is nonnegative. Therefore $c\ge 2$, and $G\setminus\{v\}$ is disconnected.
\end{proof}

 The lemma above implies that any vertex of valence $ d > g+1$ is a cut vertex. 

\begin{eg}
Given a genus $g$ rose petal graph, the number of edges emanating from the central vertex is $2g$. By adding bridges as on the right, then each vertex has valence $3$.
    
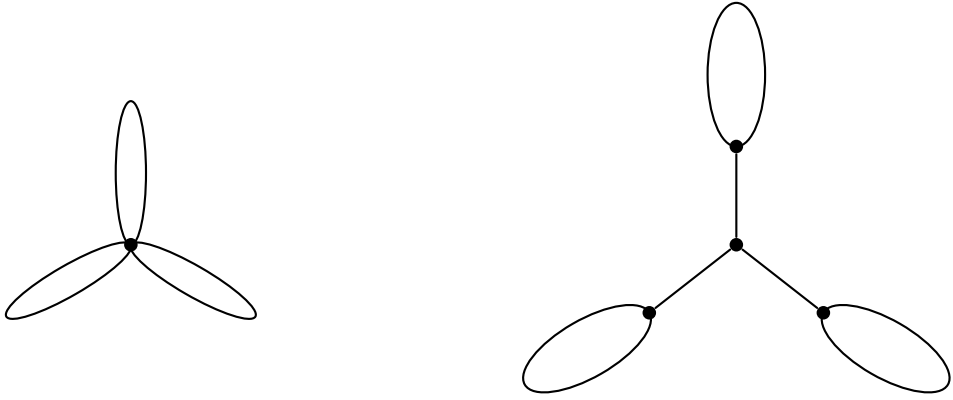
\begin{figure}[h!]
\begin{tikzpicture}[
    vertex/.style={circle,fill,inner sep=1.8pt},
    thick,
    line cap=round,
    line join=round
]


\node[vertex] (O) at (0,0) {};

\draw[rotate around={90:(0,0)}]
    (0.95,0)
    ellipse[x radius=0.95,y radius=0.2];

\draw[rotate around={210:(0,0)}]
    (0.95,0)
    ellipse[x radius=0.95,y radius=0.2];

\draw[rotate around={330:(0,0)}]
    (0.95,0)
    ellipse[x radius=0.95,y radius=0.2];


\begin{scope}[xshift=8cm]

\node[vertex] (c) at (0,0) {};
\node[vertex] (t) at (0,1.3) {};
\node[vertex] (l) at (-1.15,-0.9) {};
\node[vertex] (r) at (1.15,-0.9) {};

\draw (c)--(t);
\draw (c)--(l);
\draw (c)--(r);

\draw[rotate around={90:(t)}]
    ($(t)+(0.95,0)$)
    ellipse[x radius=0.95,y radius=0.38];

\draw[rotate around={210:(l)}]
    ($(l)+(0.95,0)$)
    ellipse[x radius=0.95,y radius=0.38];

\draw[rotate around={330:(r)}]
    ($(r)+(0.95,0)$)
    ellipse[x radius=0.95,y radius=0.38];

\end{scope}

\end{tikzpicture}
\end{figure}

\end{eg}

Given a connected graph $G$ of genus $g$, let $v$ be a vertex of valence $d>g+1$. Then one may replace $v$ by two vertices joined by a bridge so that:
\begin{enumerate}
\item the genus is unchanged;
\item contracting the bridge recovers $G$; and
\item the valence of each new vertex is strictly less than $d$.
\end{enumerate}
Indeed by Lemma \ref{disconnected}, $G\setminus\{v\}$ has at least two connected components. Partition these components into two nonempty collections. Replace $v$ by two vertices, each incident to the edges meeting one side of the partition, and join the two new vertices by a new edge $e$.

Since every path between the two sides originally passed through $v$, deleting $e$ disconnects the new graph, so $e$ is a bridge. The operation adds one vertex and one edge, so the genus is unchanged. Finally, each new vertex is incident to a proper subset of the edges formerly incident to $v$, together with the new bridge, and hence has valence strictly less than $d$. If some of the edges incident to $v$ are bridge edges, then we ignore the corresponding components in the partition above and adjoin the bridge edges to $e$ in a trivalent fashion (see Example \ref{bridgeincident}) Contracting the bridge identifies the two new vertices and recovers the original graph. One can repeat this process until the valence of each vertex of the new graph is less than or equal to $g+1$.

\begin{eg}
This example illustrates the construction above when an edge incident to $v$ is a bridge edge. 
\label{bridgeincident}
\begin{figure}[H]
\begin{minipage}{0.45\textwidth}
    \centering
    \begin{tikzpicture}[scale=1.2,
    every node/.style={circle,fill=black,inner sep=1.7pt}]

\node[label=below:$v$]  (v)  at (0,0) {};
\node[label=below:$w_0$] (w0) at (2,0) {};

\draw (v) -- (w0);

\draw (v)+(0,0.45) circle (0.45);
\draw (v)+(0,-0.45) circle (0.45);
\draw (w0)+(0.45,0) circle (0.45);
    \end{tikzpicture}
  \end{minipage}
  \hfill
\begin{minipage}{0.45\textwidth}
    \centering
    \begin{tikzpicture}[scale=1.2,
    every node/.style={circle,fill=black,inner sep=1.7pt}]

\node[label=above:$v$]  (v)  at (0,1) {};
\node[label=below:$w_0$] (w0) at (2,0) {};
\node[label=below:$v'$]  (v') at (0,-1) {};
\node[] (v'') at (0,0) {};

\draw (v'') -- (w0);
\draw (v'') --(v);
\draw (v'') -- (v');

\draw (v)+(0,0.45) circle (0.45);
\draw (v')+(0,-0.45) circle (0.45);
\draw (w0)+(0.45,0) circle (0.45);
     
    \end{tikzpicture}
  \end{minipage}
\end{figure}

\end{eg}

\begin{eg}
Consider the following genus $5$ graph. Modifying as in the above constructions requires choices, in particular how should one separate the edges emanating from the $w$ vertex. 

	
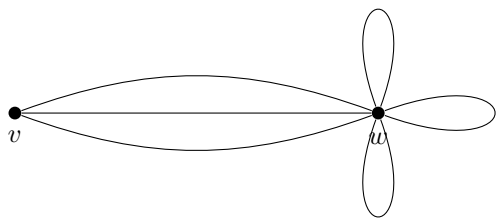
\begin{center}
	\begin{tikzpicture}[scale=1.2,
    every node/.style={circle,fill=black,inner sep=1.7pt}]

\node[label=below:$v$] (v) at (0,0) {};
\node[label=below:$w$] (w) at (4,0) {};

\draw (v) to[bend left=20] (w);
\draw (v) -- (w);
\draw (v) to[bend right=20] (w);

\draw (w) .. controls +(110:1.6cm) and +(70:1.6cm) .. (w);

\draw (w) .. controls +(20:1.8cm) and +(-20:1.8cm) .. (w);

\draw (w) .. controls +(-110:1.6cm) and +(-70:1.6cm) .. (w);

\end{tikzpicture}
\end{center}
For example, one possible separation gives
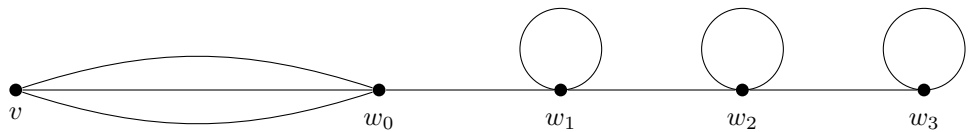
\begin{center}
	\begin{tikzpicture}[scale=1.2,
    every node/.style={circle,fill=black,inner sep=1.7pt}]

\node[label=below:$v$]  (v)  at (0,0) {};
\node[label=below:$w_0$] (w0) at (4,0) {};
\node[label=below:$w_1$] (w1) at (6,0) {};
\node[label=below:$w_2$] (w2) at (8,0) {};
\node[label=below:$w_3$] (w3) at (10,0) {};

\draw (v) to[bend left=18] (w0);
\draw (v) -- (w0);
\draw (v) to[bend right=18] (w0);

\draw (w0) -- (w1);
\draw (w1) -- (w2);
\draw (w2) -- (w3);

\draw (w1)+(0,0.45) circle (0.45);
\draw (w2)+(0,0.45) circle (0.45);
\draw (w3)+(0,0.45) circle (0.45);
\end{tikzpicture}
\end{center}
whereas a different separation would give 
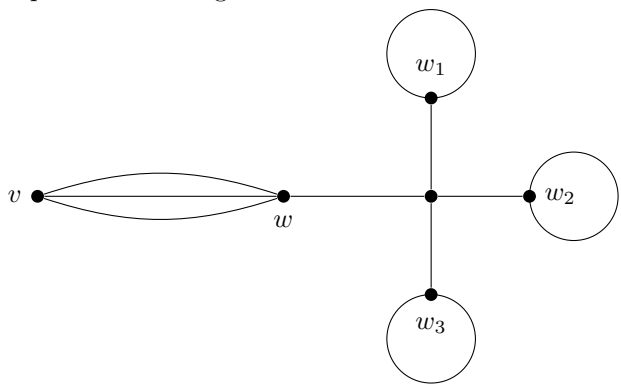
\begin{center}
\begin{tikzpicture}[scale=1.3,
    every node/.style={circle,fill=black,inner sep=1.7pt}]

\node[label=left:$v$] (v) at (-4,0) {};

\node[label=below:$w$] (w) at (-1.5,0) {};

\node (c) at (0,0) {};

\node[label=above:$w_1$]      (w1) at (0,1) {};
\node[label=right:$w_2$] (w2) at (1,0) {};
\node[label=below:$w_3$]  (w3) at (0,-1) {};

\draw (v) to[bend left=18] (w);
\draw (v) -- (w);
\draw (v) to[bend right=18] (w);

\draw (w) -- (c);

\draw (c) -- (w1);
\draw (c) -- (w2);
\draw (c) -- (w3);

\draw (w1)+(0,0.45) circle (0.45);
\draw (w2)+(0.45,0) circle (0.45);
\draw (w3)+(-0,-0.45) circle (0.45);

\end{tikzpicture}
\end{center}

\end{eg}

\begin{cor}
    Given a genus $g$ graph $G$, there exists a modification $G'$ of $G$ given by adding bridge edges such that the valence of each vertex is less than or equal to $g+1$. \end{cor}

Nevertheless, these choices of modification does not affect our constructions. Indeed, the following result is the tropical analogue of the Abel-Jacobi theorem. \begin{thm}\cite[Theorem 4.1.9]{tropicaltorelli}
    Let $\sqsubset$ and $\sqsubset'$ be  tropical curves. Then $\text{TroJac}(\sqsubset) \cong \text{TroJac}(\sqsubset')$ if and only if $\sqsubset$ and $\sqsubset'$ have the same $3$-edge connected class. 
\end{thm}

We refer the reader to \textit{loc. cit.} for the definition of $3$-edge connectiveness. We simply remark that adding bridge edges does not change the $3$-edge connected class of the graph and thus the modifications we perform above do not change the isomorphism type of the Jacobian. Moreover, by \cite[Proposition 1.12]{caparasorank}, if $\sqsubset$ and $\sqsubset$ differ by bridge edges, the isomorphism $\text{TroPic}(\sqsubset)\cong \text{TroPic}(\sqsubset')$ is rank preserving \ie if we have if $D, D'$ define the same tropical line bundle under the isomorphism, then $r_{\sqsubset}^\trop(D) = r_{\sqsubset'}^\trop(D')$. In particular, if $\sqsubset' \rightarrow \sqsubset$ is a contraction of bridge edges, then the tropical Picard groups are isomorphic and the isomorphism is rank preserving. When proving the tropical Riemann-Roch theorem in \S\ref{troprr}, we will make the following assumption. 

\begin{assumption}
Since the tropical rank is unchanged by bridge contractions, in \S\ref{troprr}, we will assume that all the bridge edges of $\sqsubset$ have been contracted. We will remark in \S\ref{globalsections} why this does not impact our computations in cohomology of line bundles on the ringed space. 
\label{contractall}
\end{assumption}

The differential of the tropical Abel map $d\alpha$ takes an integral primitive tangent vector associated to an edge $e$ of $\sqsubset$ to $$s_e := \left((\gamma_{e'})_{e'\in E(G)} \mapsto \gamma_e \right)\in N_\R = \Hom(M,\R) = \Hom(H_1(\sqsubset,\Z),\R).$$  Consider the following map: fix $v\in V(\sqsubset)$ with incident half edges $e_1,\dots,e_n$ (oriented away from $v$), and consider
$$\rho: H_1(\sqsubset,\mathbb R)\longrightarrow \mathbb R^n,\qquad \rho(\gamma)=(s_{e_1}(\gamma),\dots,s_{e_n}(\gamma)).$$
The dual of $\rho$ sends the $i$-th coordinate functional on $(\R^n)^\vee$ to $s_{e_i}$. Thus $$\dim\text{span}\{s_{e_1},\cdots, s_{e_n}\} = \rank(\rho).$$   The kernel of $\rho$ has the explicit description
$$\ker\rho=\{\gamma\in H_1(\sqsubset,\mathbb R): s_{e_i}(\gamma)=0\ \forall i\}=H_1(\sqsubset\backslash\{v\},\mathbb R),$$
and thus if $\sqsubset':= G\backslash \{v\}$ has $c$ connected components, then \begin{align*}
    \text{rank}(\rho) &= g(\sqsubset) - g(\sqsubset')= (|E(\sqsubset)|-|V(\sqsubset)| + 1) -(|E(\sqsubset')| - |V(\sqsubset')| +c) \\
    &= |E(\sqsubset)| - |V(\sqsubset)| + 1 -(|E(\sqsubset)|-n) +(|V(\sqsubset)|-1) -c\\
    &= n-c.
\end{align*} 

\begin{cor}
If $v$ is not a cut vertex, the rank of $\rho$ on $\Star(v)$ is $n-1$ where $n = |\Star(v)|.$\label{tangency}
\end{cor}

Since $\alpha$ is harmonic, $\sum_{i} s_{e_i} = 0$. Thus the $\text{im}(\rho^\vee)$ is contained in the $(n-1)$-dimensional subspace $\{\sum_i x_i = 0\}$ of $\R^n$. In particular, when $v$ is not a cut vertex, the only relation between the image of the tangent vectors in $N_\R$ is the balancing relation. 

Given the decomposition $\Jac(\sqsubset) \cong \prod_i \Jac(\sqsubset_i)$, the lattice $M=H_1(\sqsubset,\Z)$ decomposes as $M \cong \oplus M_i$ where $M_i = H_1(\sqsubset_i,\Z)$. There is an integral analogue of $\rho$, given by $\rho_\Z : H_1(G,\Z) \longrightarrow \Z^n$ whose image is isomorphic to $\check{\Lambda}_{\sqsubset_i,v}$. We make the following observation. 
\begin{cor}
	On $\sqsubset_i$, the map $M_i \longrightarrow \check{\Lambda}_{\sqsubset_i}$ is surjective. 
	\label{surjective}
\end{cor}

\begin{eg}
	In general, the map $M \rightarrow \check{\Lambda}_\sqsubset$ is not surjective. For example, if $\sqsubset$ is a bouquet of $g$ circles then $\check{\Lambda}$ at the central vertex has rank $2g$ whilst $M$ has rank $g$. 
    
\end{eg}

\begin{rmk}
    If $v$ is a not a cut vertex, then this automatically forces $n \leq g+1$. 
\end{rmk}

By Corollary \ref{tangency} and the construction above, we can find a modification of the tropical curve such that the only relation of the image of the tangent vectors under $d\alpha$ is the balancing condition.


\subsection{Maximally degenerate curves and their Jacobians}

In this section, let $K$ be a complete non-Archimedean field, and let $C$ be a smooth, projective, geometrically connected curve of genus $g$ over $K$. We say that $C$ is maximally degenerate if, after a finite extension of $K$, it has totally split degenerate reduction, i.e. the special fiber $\mathscr{C}_k$ of its minimal snc $\mathscr{C}$ model is reduced and all its irreducible components are rational. The skeleton $\Sk(C) := \Sk(\mathscr{C})$ is the dual intersection complex of  $\mathscr{C}_k$; there is a natural inclusion $\Sk(C) \hookrightarrow C^\an$ and a deformation retraction $\rho_C: C^\an \rightarrow \Sk(C)$. The skeleton $\Sk(C)$ is a tropical curve in the sense above and we will often denote the $\Sk(C)$ by $\sqsubset$. Since $C$ is maximally degenerate, $\sqsubset$ has genus $g$. 

\begin{rmk}
A smooth projective curve over a complete non-archimedean field is commonly called a \emph{Mumford curve} if it admits a Schottky uniformization. By \cite[Theorem 3.3, 4.20]{Mumford1972}, this is equivalent to being maximally degenerate. Throughout this paper, we use the term \emph{maximally degenerate} to emphasize the tropical and combinatorial aspects of the theory, although it is equivalent to the notion of a Mumford curve in our setting.
\end{rmk}

An equivalent definition of a maximally degenerate curve is that any point $x \in C^\an$ has a neighbourhood that is isomorphic to an open subset of $\proj^{1,\an}$. We thus have the following lemma. 

\begin{lem}\cite[Theorem 2.1]{put}
    Let $\mathcal{A}$ be an affinoid subdomain of $\proj^{1,\an}$. Then $\Pic(\mathcal{A}) = 0$. 
    \label{trivpic}
\end{lem}

\subsubsection{Uniformization of Jacobians}

Let $A$ be a $g$-dimensional abelian variety over $K$. The abelian variety $A$ is \textit{maximally degenerate} if, after a finite extension of $K$, the identity component of the special fiber of its Néron model is a split algebraic torus of dimension $g$. Equivalently, the abelian part of the Raynaud extension of $A$ is trivial. Thus $A \cong T/M$ where $T$ is a split $g$-dimensional torus and $M \subseteq T(K)$ is a full rank lattice. The map $T^\an \rightarrow A^\an$ is a a $K$-analytic group homomorphism whose kernel is $H_1(A^\an,\Z) \subseteq T(K)$. The tropicalization map \[\trop: T^\an \longrightarrow N_\R = \text{Hom}(M,\R)\] restricts to an isomorphism $M \xrightarrow{\sim} \trop(M)$ onto a full rank lattice $\trop(M)$ in $N_\R$. Let $\Sk(A) := N_\R/\trop(M)$; then Berkovich \cite[\S6.5]{berkovich} has shown there are natural maps $\rho_A : A^\an \rightarrow \Sk(A)$ and $\Sk(A) \rightarrow A^\an$ with a commutative diagram 
\begin{center}
	\begin{tikzcd}
		0 \ar[r] & M \ar[d] \ar[r] & T^\an \ar[r]\ar[d,"\trop"] & A^\an \ar[r] \ar[d,"\rho_A"] & 0 \\ 
		0 \ar[r]& \trop(M)\ar[r] & N_\R \ar[r] & \Sk(A) \ar[r] & 0.
	\end{tikzcd}
\end{center}
We call $\Sk(A)$ the \textit{skeleton} of $A$. The map $M \hookrightarrow T^\an \rightarrow N_\R$ defines a bilinear pairing $\langle\cdot,\cdot \rangle_\an : M \times M \rightarrow \R$ which makes the triple $(\Sk(A), M, \langle \cdot, \cdot \rangle_\an)$ a  tropical abelian variety. 
In \cite[\S6-7]{uniformisation}, it is shown that the character lattice of $T$ is $H_1(A^\an, \Z)$.
\begin{notation}
    Let $M$ denote the character lattice of $T$ and $N$ the cocharacter lattice of $T$ respectively. 
\end{notation}

By \cite[\S5]{boschlutkebohmert}, for a semistable model $\mathscr{C}$ of a curve $C$ with smooth irreducible components, there is a canonical short exact sequence 
\begin{align*}
    0 \longrightarrow \overline{T} \longrightarrow \Jac(\mathscr{C}_k) \longrightarrow \prod_i \Jac(\overline{C}_i) \longrightarrow 0
\end{align*}
where $\overline{T}$ is the split $g$ dimensional algebraic torus over $k$ and $\overline{C}_i$ are the irreducible components of $\mathscr{C}_k$. When $C$ is maximally degenerate, each $\overline{C}_i \cong \proj^1$ and thus $\Jac(\overline{C}_i) = 0$; thus the Jacobian of $C$ has totally toric reduction and thus $\Jac(C)$ is also maximally degenerate. 
The Jacobian of a maximally degenerate curve $C$ is a maximally degenerate abelian variety $\Jac(C)$ over $K$. By Proposition 4.7 in \cite{jacobian}, the homomorphism $H_1(C^\an, \Z) \rightarrow H_1(\Jac(C^\an), \Z)$ induced by the analytic Abel-Jacobi map is an isomorphism and in particular, $M \cong H_1(C^\an,\Z) = H_1(\sqsubset,\Z)$. 

We regard \(M\) as the group of cycles in the free abelian group
\(C_1(\sqsubset,\mathbb{Z})\) generated by the edges of \(\sqsubset\).
Define a pairing \(\langle \cdot,\cdot \rangle\) on
\(C_1(\sqsubset,\mathbb{Z})\) by
\begin{equation}
\tag{2.1.1}
\langle e,e' \rangle =
\begin{cases}
1 & \text{if } e=e',\\
0 & \text{otherwise}.
\end{cases}
\label{edgelength}
\end{equation}
This restricts to a symmetric, positive-definite pairing
\[
\langle \cdot,\cdot \rangle
\colon M \times M \longrightarrow \mathbb{Z},
\]
hence defines an injective homomorphism
\(\mu \colon M \to N\).

\begin{eg}
    Let $C$ be a maximally degenerate curve over $K$ and $\Jac(C)$ its Jacobian. The skeleton of $\Jac(C)$, which we denote $\text{TroJac}(\sqsubset)$, along with $M= H_1(\sqsubset,\Z)$ and the edge length pairing $\langle\cdot,\cdot\rangle$ defines a principally polarized tropical abelian variety. 
\end{eg}

The map $M \rightarrow T(K) = \text{Hom}(M,K^\times)$ is given by a pairing $[\cdot, \cdot ] : M\times M \rightarrow K^\times $. Choosing a basis for $M$, one can represent this pairing by a matrix, called the \textit{period matrix} of $\Jac(C)$. Similarly, we call the matrix representing the edge length pairing $\langle \cdot, \cdot \rangle : M \times M \rightarrow \R$ the \textit{tropical period matrix} of the principally polarized tropical abelian variety $\text{TroJac}(\sqsubset)$. By Theorem 2.9 and Remark 6.7 in \cite{jacobian}, we have that the pairing $\langle \cdot, \cdot \rangle_\an$ coincides with $\langle\cdot, \cdot, \rangle$ and $\langle \cdot ,\cdot \rangle = v_K \circ [\cdot, \cdot ] $ where $v_K: K^\times \longrightarrow \R$ is the non-archimedean valuation on $K$. 



\subsubsection{Compatibility}

The algebraic and tropical Abel-Jacobi maps are compatible in the following sense. The same statements hold true for the Abel and tropical Abel map. We will state it for the Abel map.

\begin{thm}\cite{jacobian,brandtulirsch}
Let $C$ be a smooth projective geometrically irreducible curve over $K$ and $\sqsubset$ the skeleton of $C$.  Let $\alpha_d^\an: C_d^\an \longrightarrow \Pic^{d}(C^\an) $ be the $d$-fold Abel map  $\alpha_d^\trop: \ \sqsubset_d \longrightarrow \TroPic^d(\sqsubset)$ the $d$-fold tropical Abel map as above. 
    There is a canonical isomorphism $\TroPic^d(\sqsubset) \cong \Sk(\Pic^d(C))$ such that the following square commutes
    \begin{center}
    \begin{tikzcd}
        C_d^\an \ar[r,"\alpha_d^\an"] \ar[d,"\rho_{C,d}"] & \Pic^d(C)^\an \ar[d,"\rho_{\Jac(C)}"] \\ 
        \sqsubset_d\ar[r,"\alpha_d^\trop"] & \text{TroPic}^d(\sqsubset).
    \end{tikzcd}
\end{center}
\label{jacobianskeleton}
\end{thm}

\subsubsection{Berkovich retraction and the pushforward topology}





Let $C^\an_{\rho_*}$ be the ringed space which is given set theoretically  by the analytic curve $C^\an$ equipped with the pushforward topology $\mathcal{T}_{\rho_*}$. There is a morphism of sites $\epsilon: C^\an \longrightarrow C^\an_{\mathcal{\rho_*}}$
which induces the morphism of groups 
\begin{align*}
    H^1(C^\an_{{\rho_*}}, \Oc_{C_{\rho_*}^\an}^\times ) \longrightarrow H^1(C^\an,\Oc_{C^\an}^\times). 
\end{align*}
\begin{prop}
    The morphism $H^1(C^\an_{\mathcal{T}_{\rho_*}}, \Oc_{C_{\rho_*}^\an}^\times) \longrightarrow H^1(C^\an,\Oc_{C^\an}^\times)$ is an isomorphism. 
    \label{h1iso}
\end{prop}

\begin{proof}
By Lemma \ref{trivpic}, we have $R^q\epsilon_*\G_m = 0$ for all $q>0$ and thus the Leray spectral sequence degenerates at the second page, thus giving the isomorphism. 
\end{proof}

\begin{cor}
    Consider the ringed space $(\sqsubset, \rho_*(\Oc_{C^\an}))$. Then we have \begin{align*}
        \Pic(\sqsubset,\rho_*(\Oc_\sqsubset)) = H^1(\sqsubset, (\rho_*(\Oc_{C^\an})^\times)) = \Pic(C^\an) = \Pic(C)
    \end{align*}
    \label{cor:pushforwardPic}
\end{cor}



\subsection{Essential skeleta}

Let $C$ be a maximally degenerate semistable genus $g\geq 1$ curve over $K $  and let $\Sk(C)$ denote the essential skeleton of $C$. Let $\mathscr{C}$ be the minimal snc model of $C$.  By \cite[Theorem 3.3.13]{bakernicaise}, we have $\Sk(C) = \Sk(\mathscr{C})$ and \begin{align*}
	\Sk(C) = \cup_{\omega}\Sk(C,\omega)
\end{align*} 
where the union runs over the set of non-zero regular 2-canonical forms on $C$. If we instead take the union over the non-zero global sections of $\Omega_C^1$, then $S = \cup_\omega \Sk(C,\omega)$ is equal to all the union of the closed non-bridge edges of $\Sk(C)$ by \cite[Theorem 3.4.6]{bakernicaise}. Moreover, the essential skeleton has the structure of a metric graph by \cite{bpr13}.
 By definition of the Kontsevich-Soibelman skeleton, we have that locus defining $\Sk(C,\omega)$ is given by the closure of the points $x $ such that $\wt_\omega(x) = \wt_\omega(X)$. In particular, the weight function is constant on $\Sk(C,\omega)$. By the comparison with Temkin's metric \cite[Theorem 8.3.3]{temkinmetrization}, we have that Temkin's metric is constant on $\Sk(C,\omega)$.

Recall the following construction from \cite[\S3.4]{bakernicaise}. By cohomological flatness of $\mathscr{C} \to \operatorname{Spec}(\Oc_K)$
and Grothendieck--Serre duality, the module
$H^1(\mathscr{C}, \omega_{\mathscr{C}/\Oc_K})$ is free, so that
\[
H^0(\mathscr{C}_k, \omega_{\mathscr{C}/\Oc_K}) \otimes k
\cong H^0(\mathscr{C}_k, \omega_{\mathscr{C}_k/k}).
\]
We can identify $H^0(\mathscr{C}_k, \omega_{\mathscr{C}_k/k})$
with the space of \emph{Rosenlicht differentials} on $\mathscr{C}_k$.
A Rosenlicht differential $\omega$ is, by definition, the data of
a meromorphic differential $\omega_v$ on $C_v$ for each $v \in V(G)$
such that:
\begin{enumerate}
    \item[(1)] Each $\omega_v$ has at worst logarithmic poles at
    the inverse images under $\nu$ of the singular points of
    $\mathcal{C}_k$, and it is regular everywhere else.

    \item[(2)] If $x$ is a singular point of $\mathcal{C}_k$ and
    $\nu^{-1}(x) = \{y_1, y_2\}$, then the residues of $\omega$
    at $y_1$ and $y_2$ sum to zero.
\end{enumerate}
Given $\omega \in H^0(\mathcal{C}_k, \omega_{\mathcal{C}_k/k})$
and an oriented edge $\vec{e} \in E(\Sk(C))$, let
$\operatorname{res}_{\vec{e}}(\omega)$ be the residue of $\omega$
at the point of $\widetilde{\mathcal{C}}_k$ corresponding to
$\vec{e}$. By the residue theorem, the sum
\[
\operatorname{res}(\omega)
:= \sum_{e \in E(\Sk(C))}
\operatorname{res}_{\vec{e}}(\omega)(\vec{e})
\]
belongs to $H_1(\Sk(C),k)$, so that we obtain a morphism of
$k$-vector spaces
\[
\operatorname{res}\colon
H^0(\mathscr{C}_k, \omega_{\mathscr{C}_k/k})
\longrightarrow H_1(\Sk(C),k).
\]

 By \cite[Proposition 3.5.4]{bakernicaise}, given a non-trivial simple cycle $\gamma$ in $\Sk(C)$, there exists a canonical form $\omega \in H^0(C, \Omega_C^1)$ such that $\Sk(C,\omega) = \gamma$. We can refine the result in the maximally degenerate case to give an equivalence: 

\begin{lem}
	There exists a basis $\{\omega_1,...,\omega_g\}$  of $H^0(C, \Omega_C^1)$ such that for each $i$, $\Sk(C, \omega_i)$ is a distinct simple cycle in $\Sk(C)$.
    \label{distinctcycle}
\end{lem}

\begin{proof}
	 The residue function of \cite{temkintyomkin} gives a morphism 
    \begin{align*}
        \text{res}: H^0(C,\Omega_C^1) \longrightarrow H_1(\Sk(C),K). 
    \end{align*}
By Proposition 2.4.5 in \cite{temkintyomkin}, the (graded) reduction of $\text{res}$ is equal to $\widetilde{\text{res}}$. Lemma 3.4.5 in \cite{bakernicaise} gives that $\widetilde{\text{res}}$ is an isomorphism which implies that res is an isomorphism. The result then follows from the proof of Proposition 3.5.4 in \cite{bakernicaise}. 
\end{proof}

\subsection{$K$-affine structure for abelian varieties}
\label{kaffineabelian}
From \cite{ksaffine, Karwa2024}, the $K$-affine structure $\Aff_{K,A}$ for a maximally degenerate $g$ dimensional abelian variety $A$ is given by a short exact sequence of abelian sheaves on $\Sk(A)$: 
\begin{align*}
    0 \rightarrow K^\times \rightarrow \Aff_K \rightarrow M \rightarrow 0
 \end{align*}
where $M$ is the character lattice and $M \cong \check{\Lambda}_{\Sk(A)}$, where $\check{\Lambda}_{\Sk(A)}$ is the integral cotangent bundle of $\Sk(A)$. 
The $K$-affine structure $\Aff_{K,A}$ is defined as the quotient $\frac{(\rho_A)_*(\Oc_{A^\an}^\times)}{\ker\Res_\Omega}$ where $\rho_A: A^\an \longrightarrow \Sk(A)$ is the Berkovich retraction and $\Omega \neq 0 \in H^0(A, K_A)$. The quotient is independent of the choice of scaling of the volume form $\Omega$.

 There is an alternative construction of the $K$-affine structure from which it is easier to see functoriality with the analytic structure sheaf $\Oc_{A^\an}$. Consider the constant sheaf  $K^\times \times M \subseteq (\trop)_*(\Oc_{T^\an}^\times)$ with the non-trivial action of the period lattice $L \subseteq T(K) = \Hom(M,K^\times)$ induced by the inclusion. The action has the explicit description that given $\ell \in L$, $\ell \cdot (k,m) = (\ell(m)k, m)$.
  Thus $K^\times \times M$ with this action is an equivariant sheaf on $N_\R$ and thus descends to a sheaf $(K^\times \times M)/L$ on $\Sk(A)$. The equivalence between sheaves on $\Sk(A)$ and $L$-equivariant sheaves on $N_\R$ is exact and thus preserves monomorphisms. Hence  $(K^\times \times M)_L \subseteq (\rho_A)_*(\Oc_{A^\an}^\times)$.  By \cite[\S5.3]{Karwa2024}, $(K^\times \times M)_L$ is isomorphic to $\Aff_{K,A}$.
  
  \subsubsection{Sheaves of group algebras}

\begin{defn}
    Given a sheaf of abelian groups $G$, we can define the sheaf of group algebras to be the sheaf associated to the presheaf $$U \mapsto K[G(U)] := \Z[G(U)]\otimes_\Z K.  $$ We denote the sheaf by $K[G]$. 
\end{defn}

\begin{defn}
Let $X$ be a topological space
\label{groupalgebra} and $\Fc$ a sheaf on $X$ that fits in the short exact sequence 
\begin{align*}
	0 \rightarrow K^\times \rightarrow \Fc \rightarrow \check{\Lambda} \rightarrow 0
\end{align*}
where $K^\times$ is the constant sheaf and $\check{\Lambda}$ is a sheaf of free abelian groups.  Assume that there exists a basis of $X$ such that on that basis, the short exact sequence above splits. Then 
	The presheaf \[U \longmapsto \frac{K[\Fc(U)]}{\mathcal{I}_U}\]
where $\mathcal{I}_U$ with respect to this basis is generated by local sections of the form 
\begin{align*}
 (k, \check{\lambda})\otimes 1-(1,\check{\lambda})\otimes k
\end{align*}
where $(k, \check{\lambda}) \in \Fc(U) = K^\times \times \check{\Lambda}$ with $k \in K^\times$. The sheafification is denoted $\Oc_\Fc$ and by \cite[\href{https://stacks.math.columbia.edu/tag/009N}{Tag 009N}]{stacks-project} is a sheaf on $X$. This is a sheaf of $K$-algebras and there is an embedding (of sheaves of abelian groups) $
    \Fc \rightarrow \Oc_\Fc^\times \hookrightarrow \Oc_\Fc.
$
\end{defn}

We can thus consider the group algebra $\Oc_{\Sk(A)} := 
\Oc_{\Aff_{K,A}}$ as defined in (\ref{groupalgebra}) and deduce the following result.

 \begin{thm}
	The retraction map $\rho_A: A^\an \rightarrow (\Sk(A), \Oc_{\Sk(A)})$ is a morphism of ringed spaces. 
    \label{morphismringedabelian}
\end{thm}

\begin{proof}
	By the universal property of group algebras, given $\Aff_{K,A} \rightarrow \rho_*(\Oc_{A^\an})$, there is an induced morphism $K[\Aff_{K,A}] \rightarrow \rho_*(\Oc_{A^\an})$ which naturally factors through $\Oc_{\Sk(A)}$ by construction. 
\end{proof}

Since $(\Sk(A), \Oc_{\Sk(A)})$ is a ringed space over $\spec K$, we can now discuss the $K$-valued points. 

\begin{lem}
	The set of $K$-valued points of $(\Sk(A), \Oc_{\Sk(A)})$ is equal to $(N_\R\times T(K))/L$ where the action of $L$ on $N_\R \times T(K)$ is $\ell \cdot (x,z) = (x + \trop(\ell), \ell z)$. 
	\label{Kpointsjac}
\end{lem}

As far as the authors are aware, this is the first instance where the Berkovich retraction map has been upgraded to a morphism of ringed spaces.


\subsection{Verdier Duality} 
In \cite{grossshokrieh}, the authors develop a sheaf theoretic approach to tropical homology for tropical manifolds. In this paper, we will only be interested in the case of 1-dimensional rational polyhedral spaces \ie tropical curves. We recall their results that will be of interest to us. 
\begin{notation}
	Given a topological space $X$, Let $\Z_X$ denote the constant sheaf on $X$ and $D(\Z_X)$ the bounded derived category of $\Z_X$-modules. If $\mathcal{C}^\bullet$ and $\mathcal{D}^\bullet$ are cochain complexes of $\Z_X$-modules, then $\text{Hom}_{D(\Z_X)}(\mathcal{C}^\bullet, \mathcal{D}^\bullet)$ denotes the groups of morphisms between them in $D(\Z_X)$. 
\end{notation}

\begin{defn}
    A \textit{dualizing complex} $\omega_\sqsubset^\text{Verd}$ of a tropical curve $\sqsubset$ is an element of $D(\Z_\sqsubset)$ representing the functor 
    \begin{align*}
        D(\Z_\sqsubset) &\longrightarrow D(\Z), 
        \Fc^\bullet\longrightarrow \Hom_{D(\Z)}(R\Gamma_c (\Fc^\bullet), \Z)
    \end{align*}
    where $R\Gamma_c$ is the (total) right derived functor of taking global sections with compact support, and $\Z$ is considered as a complex concentrated in degree 0. Since all tropical curves in this paper are compact, we have $R\Gamma = R\Gamma_c$. By Verdier duality, $\omega_\sqsubset^\text{Verd} = g^! \Z_*$ where $g :\ \sqsubset \rightarrow *$ is the map to a point $*$. 
\end{defn}

There is a combinatorial description of $\omega_\sqsubset^\text{Verd}$ in \cite[Remark 4.7]{grossshokrieh} given as follows. Consider the complex \begin{align*}
\oplus_{e \in {E(\sqsubset)}} \Z_e\longrightarrow \oplus_{v \in V(\sqsubset)} \Z_v
\end{align*}
placed in degrees $[-1, 0]$. To define the differentials, one needs to pick orientations on all edges on $\sqsubset$. With respect to a fixed orientation, the differential is $0$ between the components $\Z_e$ and $\Z_v$ if $v \not \subseteq e$ and multiplication by $\epsilon_{e/v}$ where $\epsilon_{e/v} = 1$ if the edge $e$ is outgoing and $-1$ if the edge is incoming.


\begin{thm}\cite[Corollary 6.9]{grossshokrieh}
    There is a natural isomorphism 
    \begin{align*}
      \omega_\sqsubset^\text{Verd} \cong  \check{\Lambda}_\sqsubset[1] 
    \end{align*}
    where $\check{\Lambda}[1] \in D(\Z_\sqsubset)$. 
\end{thm}



%
 
\subsection{Ringed Spaces}
We collect several general results and constructions for ringed spaces that will be used later. 
Let $\textbf{RS}$ denote the category of ringed spaces and  $(X,\Oc_X) \in \textbf{RS}$. We review the construction of the categorical quotient in \textbf{RS} from \cite[Chapter 6]{dolgachev}.  One can form the product $X^n$ in $\textbf{RS}$ given by the underlying Cartesian product of the topological space and the tensor product of the sheaf of rings. The symmetric group $S_n$ acts on $X^n$ and one can take the categorical quotient $X_n  := X^n/S_n$ in $\textbf{RS}$ which is simply the ringed space whose topological space is the quotient $|X^n|/G$ and whose sheaf of rings is $q_*(\Oc_{X^n})^{S_n}$ where $q: |X^n| \rightarrow |X^n|/S_n$.

\begin{lem}
    Let $f : X \rightarrow Y$ be a morphism in \textbf{RS}. This induces a morphism $f_n : X_n \rightarrow Y_n$ in \textbf{RS}.
    \label{symmetricmorphism}
\end{lem}

\begin{proof}
    This follows from the universal property of categorical quotients.
\end{proof}



\subsubsection{Duality for ringed spaces}

For later use, we recall Spaltenstein’s duality theorem for unbounded complexes of sheaves on ringed spaces \cite[Theorem B]{spaltenstein}. As we are working with tropical curves, we satisfy the assumptions of the theorem. 

\begin{thm}
    Let the ringed spaces \(X\), \(Y\) be locally compact and locally finite dimensional with $f: X \rightarrow Y$ a morphism of ringed spaces. Then the functor \(Rf_!\) has an adjoint \(f^!\), and for
\(\mathcal A^\bullet \in \mathcal D(X)\),
\(\mathcal B^\bullet \in \mathcal D(Y)\), we have
\begin{enumerate}
	\item There is an equivalence\[
R\mathcal{H}om^\bullet(Rf_!\mathcal A^\bullet,\mathcal B^\bullet)
\simeq
Rf_*
R\mathcal{H}om^\bullet(\mathcal A^\bullet,f^!\mathcal B^\bullet).
\]
\item for $\mathcal{C}^\bullet \in \mathcal{D}(Y)$, there is a natural isomorphism \[
Rf_!(\mathcal{A}^\bullet \otimes^L_{\Oc_X} Lf^*\mathcal{B}^\bullet) \cong Rf_! \mathcal{A}^\bullet \otimes^L_{\Oc_Y}\mathcal{B}^\bullet
\]
\end{enumerate}

\label{grothendieck duality}

\end{thm}

\section{$K$-Tropical curves}

In this section, we will build the theory of $K$-tropical curves.  We will give the construction of the $K$-affine structure discussed in the introduction of this paper in \S\ref{kaffinestruc}. In \S\ref{ktropsheaf}, we give both constructions. In \S\ref{analyticmod},  we modify the structure sheaf arising from a $K$-affine structure by a suitable completion to ensure the cohomology groups of line bundles are non-trivial and to ensure that the dualizing sheaf is concentrated in degree one.

In \S\ref{ktropdivisor}, we introduce the notion of a $K$-tropical Weil divisor and define linear equivalence via the Abel map.  We end the section by proving a $K$-tropical analogue of Theorem 1.1 of \cite{jacobian}, that is the $K$-tropical principal divisors surject onto the tropical principal divisors of $\underline{\sqsubset}$. As a corollary, we give a $K$-tropical analogue of the Jacobi inversion theorem.  In \ref{divisorlinebundle}, we explicitly describe effective Cartier divisors on the modified ringed space in terms of Weil divisors whose support is at 2-valent points of $\underline{\sqsubset}$. This gives a correspondence between line bundles and divisors. In \S\ref{globalsections}, we prove that the cohomology of line bundles is finite dimensional, introduce the rank of a linear system of a divisor and relate it to the cohomology of the corresponding line bundle. Furthermore, we prove a $K$-tropical specialization lemma and give a GAGA-type result between the Berkovich analytification of a curve and the corresponding $K$-tropical curve.

\subsection{$K$-affine structure}
\label{kaffinestruc}

In this section, we will define the $K$-affine structure for a tropical curve. The construction was alluded to in the introduction but here we make it more explicit. When we consider $\underline{\sqsubset}$ as arising from the essential skeleton of an algebraic curve, we will directly define a ringed space structure on $\underline{\sqsubset}$. We will refer to the former situation as the \emph{tropical case} and the latter as the \emph{algebraic case}. 
\subsubsection{Construction for a tropical curve}
Let $\underline{\sqsubset}$ be a tropical curve and $\text{TroJac}(\sqsubset)$ the tropical Jacobian of $\underline{\sqsubset}$. Consider the section $s$ of $\val: K^\times = \C((t^\R))^\times \longrightarrow \R$ given by $a \mapsto t^a$. This section lifts the tropical period matrix of $\text{TroJac}(\sqsubset)$ to a period matrix over $K$ and thus defines a maximally degenerate abelian variety over $K$ whose skeleton is exactly $\text{TroJac}(\sqsubset)$. Following \S\ref{kaffineabelian}, we can thus construct the $K$-affine structure $\Aff_{K,\text{TroJac}(\sqsubset)}$ on $\text{TroJac}(\sqsubset)$ and form the ringed space which we denote $\Jac(\sqsubset)$. We will abuse notation and write $\Jac(C)$ for the maximally degenerate abelian variety associated to the lift of the period matrix. 

\begin{lem}
   If the tropical Jacobian $\TroJac(\sqsubset)$ decomposes as $\TroJac(\sqsubset_1)\times \dots \times \TroJac(\sqsubset_r)$ where $\underline{\sqsubset}_i$ are maximal 2-connected subgraphs, then $\Jac(C)$ is not in the Torelli locus $\mathcal{T}_g \subseteq \mathcal{A}_g$ of principally polarized $g$ dimensional abelian varieties. 
   \label{notTorelli}
\end{lem}

Before we prove the lemma, we recall the following classical result. 

\begin{lem}
    Let $C$ be a smooth projective irreducible curve of genus $g$ over $K$. Then $\Jac(C)$ is indecomposable as a principally polarized abelian variety. 
    \label{indecompose}
\end{lem}

\begin{proof}By Riemann’s theorem,
$\Theta_C=W_{g-1}(C)$
is irreducible, since it is the image of the irreducible variety $C^{(g-1)}$. On the other hand, a principally polarized abelian variety is decomposable if and only if its theta divisor is reducible. Hence $\Jac(C)$ is indecomposable. 
\end{proof}

\begin{proof}(of Lemma \ref{notTorelli})
If $\TroJac(\sqsubset)$ is decomposable, then the tropical period matrix is block diagonal. Lifting the tropical period matrix via the section tells us that the lifted period matrix is also block diagonal (over $K^\times)$ and thus the abelian variety decomposes into a product of lower dimensional principally polarized abelian varieties. 
\end{proof}

The tropical curve comes equipped with a sheaf $\Aff_{\Z,\sqsubset}$ of integral affine functions which fits in the short exact sequence 
\begin{align*}
    0 \longrightarrow \R \longrightarrow \Aff_{\Z,\sqsubset} \longrightarrow \check{\Lambda}_\sqsubset \longrightarrow 0 .
\end{align*}
We give the following construction that modifies the integral affine structure on $\underline{\sqsubset}$. This will not change the Tropical Picard group, the tropical divisor theory and thus for the proof of the Tropical Riemann-Roch theorem, it is sufficient to work with this modification. This modification will prove useful when computing the dualizing complex in \S\ref{dualizingsheaf}. \begin{construct}
Suppose $v$ is a vertex of $\underline{\sqsubset}$ of valence strictly greater than $g+1$. By Lemma \ref{disconnected}, $v$ is a cut vertex and thus by the modification considered in $\S \ref{tropmod}$, one can add bridge edges to reduce the valence of $v$. In particular, the tropical Jacobian $\TroJac(\sqsubset)$ decomposes. Let $\sqsubset_1, \cdots, \sqsubset_r$ denote  the maximal 2-connected subgraphs of $\underline{\sqsubset}$ and $\check{\Lambda}_{\sqsubset_i}$ denote the cotangent bundle of $\sqsubset_i$; it is naturally a primitive sublattice of $\check{\Lambda}_\sqsubset$. Furthermore, there is a natural inclusion $\oplus_i \check{\Lambda}_{\sqsubset_i} \subseteq \check{\Lambda}_\sqsubset$. By contravariance of $\text{Ext}^1$ in the second entry, there is an induced map \begin{align*}
    \text{Ext}^1(\R, \check{\Lambda}_\sqsubset) \longrightarrow \text{Ext}^1(\R, \oplus_i\check{\Lambda}_{\sqsubset_i}). 
\end{align*}
The image of $\Aff_{\Z,\sqsubset}$ under this map defines a modifed integral affine structure on $\sqsubset$. We will work with this modified integral affine structure in this paper; to avoid complicated notation, we will continue to denote this modified integral affine structure by $\Aff_{\Z,\sqsubset}$ and $\check{\Lambda}_\sqsubset$ will now denote $\oplus_i\check{\Lambda}_{\sqsubset_i}$. If every vertex of $\underline{\sqsubset}$ has valence less than or equal to $g+1$, then we do not do anything. 

Since the tropical Jacobian decomposes according to the $\sqsubset_i$, by construction, the tropical Picard group of this modified tropical curve is unchanged and further the induced map on cotangent lattices given by the Abel map factors through the modified cotangent lattice. Since we have also not changed the underlying graph, the notion of divisor and linear equivalence remains the same. 
\label{modifedaffinestructure}
\end{construct}

We then define $\Aff_{K,\sqsubset}:= \frac{K^\times \times \Aff_{\Z,\sqsubset}}{\{(t^n,-n): n \in \R\}}$. This then fits in the short exact sequence 
\begin{align*}
    0 \longrightarrow K^\times \longrightarrow \Aff_{K,\sqsubset} \longrightarrow \check{\Lambda}_\sqsubset \longrightarrow 0
\end{align*}
which is compatible with the short exact sequence above. 
The same construction can be applied to the sheaf of integral affine functions $\Aff_{\Z,\text{TroJac}(\sqsubset)}$ which would give us the $K$-affine structure $\Aff_{K, \text{TroJac}(\sqsubset)}$. By functoriality of pushout, the map $\Aff_{\Z,\text{TroJac}(\sqsubset)} \longrightarrow \alpha_*\Aff_{\Z, \sqsubset}$ gives a map \begin{align*}
    \Aff_{K,\text{TroJac}(\sqsubset)} \longrightarrow \alpha_*\Aff_{K, \sqsubset}.
\end{align*} 
 
It fits in the following diagram
 	\begin{center}
 		\begin{tikzcd}
 			0 \ar[r] & K^\times \ar[r]\ar[d,"\text{id}"] &\Aff_{K,\TroJac(\sqsubset)}\ar[r]\ar[d] & \ar[r] \ar[d] M & 0\\
 			0 \ar[r] & K^\times\ar[r] & \alpha_*^\trop\Aff_{K,\sqsubset} \ar[r] & \alpha_*^\trop\check{\Lambda}_\sqsubset \ar[r] & 0.
 		\end{tikzcd}
 	\end{center}

Taking long exact sequences and given the isomorphisms $$H^0(K^\times) \cong K^\times \cong  H^0(\Aff_{K,\TroJac(\sqsubset)})\cong H^0(\Aff_{K,\sqsubset})$$ gives the commutative diagram \begin{center}
	\begin{tikzcd}
		0 \ar[r] &H^0(M) \ar[r]\ar[d] & H^1(K^\times) \ar[d,"\text{id}"]  &\\
		0 \ar[r] & H^0(\check{\Lambda}_\sqsubset)  \ar[r]& H^1(K^\times) \ar[r] &H^1(\Aff_{K,\sqsubset}) \ar[r] & H^1(\check{\Lambda}_\sqsubset) \ar[r] & 0.
	\end{tikzcd}
\end{center}
The vertical arrow $H^0(M) \rightarrow H^0(\check{\Lambda}_\sqsubset)$ is the canonical identification of the integral cotangent bundle of the tropical Jacobian with the harmonic 1-forms on $\underline{\sqsubset}$.  In particular this proves: 
\begin{prop}
	Given $\underline{\sqsubset}$ and $\Jac(C)$ as above, there is an isomorphism $H^1(\Aff_{K,\sqsubset}) \cong \Z \times \Jac(C)(K)$. 
	\label{pic}
\end{prop}

\subsubsection{The algebraic case}

In the algebraic case, the ringed space structure does not arise from a $K$-affine structure. In fact, in Example \ref{ex:noKaffinealg}, we show that such a $K$-affine structure cannot exist that is also compatible with the $K$-affine structure associated to the Jacobian of $C$.  In future work, the authors will give an alternative construction such that given any maximally degenerate abelian variety $A$, there exists a $K$-tropical curve whose Jacobian is exactly $A$. In particular, the Torelli locus for $K$-tropical curves is all maximally degenerate principally polarized abelian varieties.

\begin{eg}
\label{ex:noKaffinealg}
    Let $C/K$ be a smooth, projective,
geometrically connected Mumford curve of genus two and suppose that its essential skeleton
$\underline{\sqsubset}$ is a dumbbell: two loops $L_1,L_2$, based at vertices $v_1,v_2$,
joined by a bridge $b$. 

Choose oriented generators $\ell_1,\ell_2$ corresponding to the two
loops and let $u_1,u_2\in M$ denote the corresponding basis of characters in $M$.  Let $\Jac(C)^\an \cong \G_m^{2,\an}/L$ and let $Q=(Q_{ij})_{1\leq i,j\leq2}$ be the period matrix of $C$ with respect to this basis.

Let $\Aff_{K,\TroJac(\sqsubset)}$ denote the $K$-affine structure sheaf on $\TroJac(\sqsubset)$ constructed in \S\ref{kaffineabelian}.
We then claim that there is no extension of sheaves of abelian groups on $\sqsubset$
\[
  1\longrightarrow{K^\times}
  \longrightarrow\Aff_{K,\sqsubset}
  \longrightarrow{}\check\Lambda_\sqsubset
  \longrightarrow0
\]
admitting a morphism
\[
\Phi: \Aff_{K,\TroJac(\sqsubset)}\longrightarrow\alpha_*^\trop\Aff_{K,\sqsubset}
\]
such that there is a
commutative diagram with exact rows
\[
\begin{tikzcd}[column sep=small]
1 \arrow[r] & {K^\times} \arrow[r] \arrow[d,equal]
  & \Aff_{K,\TroJac(\sqsubset)} \arrow[r] \arrow[d]
  & \check\Lambda_{\TroJac(\sqsubset)} \arrow[r] \arrow[d,"h"] & 0 \\
1 \arrow[r] & {K^\times} \arrow[r]
  & \alpha_*^\trop\Aff_{K,\sqsubset} \arrow[r]
  & \alpha_*^\trop\check\Lambda_\sqsubset \arrow[r] & 0.
\end{tikzcd}
\]
Indeed, the harmonic form corresponding to $u_2$ is $h(u_2)$. Let $U$ be the open neighbourhood of
$L_1$ consisting of the entire left loop together with a short initial
segment of the bridge, excluding $v_2$. Then $
  h(u_2)|_U=0.$
In particular, the slope vanishes at $v_1$ in all three outgoing
directions.

Consider the sheaf $\mathcal P$ of lifts of $u_2$ to
$(\alpha^\trop)^{-1}\Aff_{K,\TroJac(\sqsubset)}|_U$. It is a torsor under
${K^\times}$, with monodromy $Q_{21}$ around $L_1$.
If $\Phi$ existed, then $h(u_2)|_U=0$ would imply that $\Phi$ sends every
 section of $\mathcal P(U)$ to an element of
${K^\times}$. Because $\Phi$ is the identity on constants, this
map of torsors is equivariant and trivializes $\mathcal P$.
Consequently the monodromy must be trivial:
\begin{equation}\label{eq:cross-period}
  Q_{21}=Q_{12}=1.
\end{equation}
This would give a decomposition of $\Jac(C)$ but this is not possible by Lemma \ref{indecompose}. 
\end{eg}

\subsection{$K$-tropical structure sheaf}
\label{ktropsheaf}


In the tropical case, we define the $K$-tropical structure sheaf as follows. 

\begin{defn}
	The $K$-tropical structure sheaf $\Oc_\sqsubset$ is defined as $\Oc_{\Aff_{K}}$ following (\ref{groupalgebra}). We then define the ringed space $\sqsubset:= (\underline{\sqsubset},\Oc_\sqsubset)$ and refer to it as a $K$-tropical curve. 
\end{defn}

The following proposition follows readily from the definition of $\Oc_\sqsubset$. 

\begin{prop}
    We have a canonical isomorphism $\Aff_{K,\sqsubset}\cong \Oc_{\sqsubset}^{\times} $ as sheaves of abelian groups. 
    \label{units}
\end{prop}
In \S\ref{analyticmod}, we will modify this sheaf by completing it. This is not necessary in the algebraic case where we define the $K$-tropical structure sheaf as follows: 

\begin{defn}
    The $K$-tropical structure sheaf in the algebraic case is defined as $\Oc_\sqsubset^\an := \rho_*(\Oc_{C^\an})$ where $\rho: C^\an \longrightarrow \underline{\sqsubset}$ is the retraction map. 
\end{defn}

\begin{rmk}
   In the tropical case, we have assumed that we contract all the bridge edges. This will also be beneficial for when we compute the dualizing complex of $(\underline{\sqsubset},\Oc_\sqsubset^\an)$ in the algebraic case. We will denote the contraction by $c : \underline{\sqsubset}\longrightarrow \underline{\sqsubset}_\text{nb}$ and consider the ringed space $(\underline{\sqsubset}_\text{nb}, c_*\Oc_{\sqsubset^\an})$. 
\end{rmk}

The pair $\sqsubset= (\underline{\sqsubset},\Oc_\sqsubset)$  is a ringed space and we can thus discuss modules over the structure sheaf $\Oc_\sqsubset$. In this paper, we will be solely interested in invertible sheaves on $\sqsubset$. We denote the Picard group of $\sqsubset$ by $\Pic(\sqsubset)$, which is not to be confused with the tropical Picard group $\text{TroPic}(\sqsubset)$. 

\begin{cor}
In the tropical case, the Picard group $
\Pic(\sqsubset)$ is isomorphic to $\Z\times \Jac(C)(K)$. 
\label{pic2}
\end{cor}

\begin{proof}
	This follows immediately from Propositions \ref{pic} and \ref{units}. 
\end{proof}

\begin{lem}
    In the algebraic case, the Picard group $\Pic(\underline{\sqsubset},\Oc_\sqsubset^\an)$ is isomorphic to $\Pic(C)$. The same holds true for the ringed space $(\underline{\sqsubset}_\text{nb}, c_*\Oc_{\sqsubset}^\an)$. 
\end{lem}

\begin{proof}
    This follows immediately from Corollary \ref{cor:pushforwardPic}.
\end{proof}

\subsubsection{Compatibility with Berkovich analytic category}
Given the ringed space structure in the algebraic case, we have the following theorem. 
\label{compatible}

\begin{thm}\label{thm:modified-ringed-space-morphism}

 Let $C^\an$ be the Berkovich analytification of a Mumford curve and $\underline{\sqsubset}$ the essential skeleton of $C$. Then $(\underline{\sqsubset},\Oc_{\sqsubset}^\an)$ is a ringed space such that 
 
 \begin{enumerate}
 	\item     The retraction map $\rho$ induces a morphism
of ringed spaces
\[
    (C^{\an},\Oc_{C^{\an}})
    \longrightarrow
    (\underline{\sqsubset},{\Oc}_{\sqsubset}^\an).
\] 
 	\item The pullback $\rho_C^*$ induces an isomorphism $\Pic(C^\an) \cong \Pic(C) \longrightarrow \Pic(\underline{\sqsubset},\Oc_\sqsubset^\an)$.
 	\item The $K$-tropical Jacobian $\Jac(\sqsubset):= (\TroJac(\sqsubset),\Oc_{\Jac(\sqsubset)})$ is a ringed space such that the tropical Abel-Jacobi map $\alpha : \sqsubset \rightarrow \Jac(\sqsubset)$ is a morphism of ringed spaces and there is a commutative diagram of ringed spaces 
        \begin{center}
            \begin{tikzcd}
                C^\an \ar[r, "\alpha^\an"] \ar[d,"\rho_C"] & \text{Jac}(C^\an) \ar[d, "\rho_{\Jac(C)}"] \\ 
                \sqsubset \ar[r, "\alpha^\trop"] &\text{Jac}(\sqsubset)
            \end{tikzcd}
        \end{center}
        where $\alpha^\an$ is the analytic Abel-Jacobi map. 
       
 \end{enumerate}
\label{algebraicKtropical}

\end{thm}

	%
%

\subsection{The analytic modification of the structure sheaf}

\label{analyticmod}

The current construction of the structure sheaf $\Oc_\sqsubset$ in the tropical case is not suitable for our purposes as illustrated by the following example. 

\begin{eg}
    Let $E = \G_m/q^\Z$ be a Tate curve. Then $\underline{\sqsubset}$ is a circle of length $-\val(q)$. The structure sheaf $\Oc_\sqsubset$ is locally isomorphic to $A = K[z^{\pm1}]$ with transition function given by $z \mapsto qz$. Consider the degree $d$ line bundle $\Lc_d$ whose transition function is given by $z^d$ for $d\geq 0$. Then \begin{align*}
        H^0(\sqsubset, \Lc_d) = \ker (A  \rightarrow A )
    \end{align*}
    where the map is defined as $z \mapsto qz^{d+1} - z$. The kernel of this map is always zero and thus there are no global sections of $\Lc_d$. This is contrary to the algebraic analogue where every effective divisor has a non-zero global section. 
\end{eg}

\begin{rmk}
    The key reason that the previous example of a positive degree line bundle on a Tate curve fails is that the line bundles do not have theta functions as global sections.
\end{rmk}

To fix this, we replace our of sheaf of rings which are locally polynomial algebras with their completions. Let $e$ be an edge and $\mathring{e}$ denote its interior. Let $A_e := \Oc_\sqsubset|_{\mathring{e}}$ and consider \(S\), a compact interval contained in $\mathring{e}$.  Let \(A_S := \Oc_{\sqsubset}|_S\). The \(K\)-affine structure associates to each unit \(u\in A_S^\times\) an affine function \(\ell_u:S\to\mathbb R\), with
\[
\ell_{uv}=\ell_u+\ell_v,
\qquad
\ell_c=-\log|c|\quad(c\in K^\times).
\]
For \(f\ne0\), write \(f=\sum_j u_j\) as a sum of monomial terms whose classes in \(A_S^\times/K^\times\) are distinct. Define
\[
|f|_x:=\max_j e^{-\ell_{u_j}(x)}
\qquad(x\in S),
\]
and set \(|0|_x=0\). Then
\[
\|f\|_S:=\sup_{x\in S}|f|_x
\]
defines a non-archimedean norm on \(A_S\).  We denote its completion with respect to this norm by \(\widehat A_S\). For compact intervals \(S'\subseteq S\), the restriction map satisfies
\[
\|f|_{S'}\|_{S'}\leq\|f\|_S
\]
and therefore extends uniquely to a continuous homomorphism \(\widehat A_S\to\widehat A_{S'}\). Therefore 
    the norms $
\{\|\cdot\|_{S}\}_{S\Subset\mathring{e}}$ form a directed family and define a locally convex
non-Archimedean topology on $A_e$.


\begin{defn}
    The algebra associated to the open edge $\mathring{e}$ is the 
completion of $A_e$ with respect to the family of norms arising
from all compact intervals:
\[
\Oc_\sqsubset^\an(\mathring{e})= A_e^{\an}(\mathring{e})
:=
\widehat{A_e}^{\,
\{\|\cdot\|_{S}\}_{S\Subset \mathring{e}}}.
\]
For $S'\subseteq S'$, the inequality $
\|\cdot\|_{S'}\leq \|\cdot\|_{S}$
induces a continuous restriction homomorphism
\[
A_e^{\an}(S)\longrightarrow A_e^{\an}(S').
\]
The algebras $A_e^{\an}(S)$ therefore form an inverse system, and the algebra on $\mathring{e}$ can equivalently be written as
\begin{align}
A_e^{\an}(\mathring{e})
\simeq
\varprojlim_{S\Subset \mathring{e}}A_e^{\an}(S)
=
\varprojlim_{S\Subset \mathring{e}}
\widehat{A_e}^{\,\|\cdot\|_{S}}.
\end{align}
\label{completeedge}
\end{defn}

\begin{rmk}
    The construction also works for the open star of a 2-valent vertex where if the open star of $v$ is given by say $(0,r_1] \cup [0,r_2)$, then we take the completion with respect to all closed intervals in $(0,r_1+r_2)$.
\end{rmk}

We now define the completion of $\Oc_{\sqsubset,v}$ by imposing convergence along each branch emanating
from $v$. Let $S_v =  \cup_{i=1}^d [v,x_i]$ be a closed star neighbourhood of $v$ and write $A_v :=\Oc_{\sqsubset}|_{S_v}$. Write $e_i = [v,x_i]$ for the branches of $S_v$ and $A_i := A_{e_i}$. 
Let $C_v$ be an affinoid domain whose skeleton is precisely $S_v$. By \cite{Mehmeti2019PatchingOB}, for an affinoid domain, the topological boundary of $C_v$ coincides with the Shilov boundary. Thus the boundary is $\{x_1,\cdots, x_d\}$. 
The restriction map to ${e}_i \cap S_v$  gives the
$K$-algebra homomorphism $
\text{res}_{v,e_i}:A_v\longrightarrow A_i^\an.$ Then we define
\[
\|f\|_{S_v}
:=
\max_{1\leq i\leq d}\|\text{res}_{v,e_i}(f)\|_{A_i},
\qquad f\in A_v.
\]
where $||\cdot||_{A_i}$ is the norm on $A_i^\an$. 
Since each $\|\cdot\|_{A_i}$ is a non-Archimedean submultiplicative
norm, $\|\cdot\|_{S_v}$ is a non-Archimedean submultiplicative
seminorm on $A_v$. In particular,
\[
\|f+g\|_{S_v}
\leq
\max\{\|f\|_{S_v},\|g\|_{S_v}\}
, \quad 
\|fg\|_{S_v}
\leq
\|f\|_{S_v}\|g\|_{S_v}.
\]

\begin{lem}
     Let $\mathcal{T}_v := \ker \left( A_v \longrightarrow \oplus_i A_{i}\right)$. The seminorm $||\cdot||_{S_v}$ is a norm on $A_v /\mathcal{T}_v$. 
\end{lem}

\begin{proof}
    Since $||\cdot||_{A_i}$ is a norm, $||f||_{S_v} = 0$ if and only if $\text{res}_{v,e_i}(f) =0$ for every $i$ \ie $f \in \mathcal{T}_v$. 

\end{proof}

We now pass from closed stars to open star neighbourhoods. Let $U_v$ be
an open star neighbourhood of $v$. Let $S' \Subset U_v$ be a compact star.   For each such $S'$ the construction
above gives a seminorm $\|\cdot\|_{S'}$
on $A_v$ such that if $S'\subseteq S''$, then
$\|f\|_{S'}\leq \|f\|_{S''}$
or every $f\in A_v$. Therefore we have, 
\begin{lem}
    The seminorms $
\{\|\cdot\|_{S'}\}_{S'\Subset U_v}$
therefore form a directed family and define a locally convex
non-Archimedean topology on $A_v$.
\end{lem}

We can thus make the following definition. 
\begin{defn}
\label{completestar}
The  algebra associated to the open star $U_v$ is the 
completion of $A_v$ with respect to the family of seminorms arising
from all compact substars:
\[
\Oc_\sqsubset^\an(U_v)= A_v^{\an}(U_v)
:=
\widehat{A_v/\mathcal{T}_v}^{\,
\{\|\cdot\|_{S'}\}_{S'\Subset U_v}}.
\]
For $S'\subseteq S''$, the inequality $
\|\cdot\|_{S'}\leq \|\cdot\|_{S''}$
induces a continuous restriction homomorphism
\[
A_v^{\an}(S'')\longrightarrow A_v^{\an}(S').
\]
The algebras $A_v^{\an}(S')$ therefore form an inverse system, and the algebra on $U_v$ can equivalently be written as
\begin{align}
A_v^{\an}(U_v)
\simeq
\varprojlim_{S'\Subset U_v}A_v^{\an}(S')
=
\varprojlim_{S'\Subset U_v}
\widehat{A_v/\mathcal{T}_v}^{\,\|\cdot\|_{S'}}.
\end{align}\end{defn}
By definition of the norm $||\cdot||_{S_v}$, for every $i$, we have $
\|\text{res}_{v,e_i}(f)\|_{A_i}\le \|f\|_{S_v}$. Thus \(\text{res}_{v,e_i}\) is contractive and  sends the kernel of a vertex seminorm into the kernel of a branch seminorm. It therefore descends to the separated completions and extends uniquely to a contractive homomorphism between the respective completions. 
These local models and restriction maps give rise to a sheaf on $\sqsubset$.


\begin{defn}
    Let $\Oc_\sqsubset^\an$ denote the sheaf arising from the completion of $\Oc_\sqsubset$ described above. We define $\sqsubset^\an := (\underline{\sqsubset},\Oc_\sqsubset^\an)$ and shall also refer to it as a $K$-tropical curve for simplicity. 
\end{defn}



We next describe the kernel $\mathcal{T}_v$ as it will be used in \S\ref{dualizingsheaf} when computing the dualizing complex. We first introduce the notion of elliptic $m$-fold singularity. 

\begin{definition}
 For $m\geq 3$, an elliptic $m$-fold singularity over $K$ is a ring isomorphic to

\[
R_m \simeq
\begin{cases}
K[x_1,x_2]/(x_1^2x_2-x_1x_2^2),
& m=3 \quad \text{(planar triple point)},\\
K[x_1,\ldots,x_{m-1}]/I_m,
& m\geq 4 \quad
\text{(\(m\) general lines through the origin in \(\mathbb A^{m-1}\))},
\end{cases}
\]
where
\[
I_m
:=
\left(
x_hx_i-x_hx_j
\;\middle|\;
i,j,h\in\{1,\ldots,m-1\}\text{ distinct}
\right).
\]
The ring $R_m$ is a one-dimensional reduced Gorenstein 
ring with $m$ branches by \cite[Appendix A]{smyth}. 
\end{definition}

By a coordinate change $x_i \mapsto x_i-1$, the localization of $R_m$ at the ideal $(x_1-1,\cdots, x_{m-1}-1)$ is isomorphic to 
\[R_m^\times =
\begin{cases}
K[x_1^{\pm1},x_2^{\pm1}]/((x_1-1)(x_2-1)(x_1-x_2)),
& m=3 \\
K[x_1^{\pm1},\ldots,x_{m-1}^{\pm1}]/I_m^\times,
& m\geq 4 
\end{cases}
\]
where
\[
I_m^\times
:=
\left(
(x_h-1)(x_i-x_j)
\;\middle|\;
i,j,h\in\{1,\ldots,m-1\}
\text{ are pairwise distinct}
\right)
\]
Note the ring $R_m^\times$ is also Gorenstein.

\begin{lem}
\label{kernel}
    Let $d$ denote the valence of the vertex $v$ and take the following basis $$u_1 := e_1 - e_{d}, u_2:= e_2 - e_d, \cdots, u_{d-1}:= e_{d-1}-e_d$$ for $\check{\Lambda}_{\sqsubset,v}$. Then the kernel of $A_v \longrightarrow \oplus_i A_i$ is generated by 
    \begin{align*}
        \bigcap_{i-1}^{d-1}\left((\chi^{u_j}-1:j\neq i)\cap \left(\chi^{u_i}-\chi^{u_j}: 1 \leq i,j\leq d-1\right)\right)
    \end{align*} 
   and $A_v/\mathcal{T}_v$ is isomorphic to $R^\times_d$. In particular, $A_v/\mathcal{T}_v$ is Gorenstein. 
\end{lem}

\begin{proof}
This follows from the fact that $I_m$, and hence $I_m^\times$, is equivalently described as the kernel of restricting to every branch of the elliptic $m$-fold singularity as explained in the discussion in \cite{smyth} after Definition 2.1 and thus agrees with the kernel of $A_v \longrightarrow \oplus_iA_i$ too. 
\end{proof}

\begin{prop}
Let $S_v$ be a closed star of the vertex $v$ and assume that all the branches of $S_v$ have the same length. 
The completion $\widehat{A}_v^{\|\cdot\|_{S_v}}$ of $A_v$ with respect to $||\cdot||_{S_v}$ is
Gorenstein and the natural map
$A_v\to \widehat{A}_v^{\|\cdot\|_{S_v}}$ is flat.
\label{gorenstein}
\end{prop}

\begin{proof}
We claim there is a canonical isomorphism of topological
$K$-algebras
\begin{align}
    \widehat A_v^{\,\|\cdot\|_{S_v}}
\cong
\mathscr S_{r,R}/I_m^\times\mathscr S_{r,R},
\label{ellipticaffinoid}
\end{align}
where $
\mathscr{S} 
=
K\langle
R^{-1}x_1,\ldots,R^{-1}x_{m-1},\,
rx_1^{-1},\ldots,rx_{m-1}^{-1}
\rangle
$
with $0 < r\leq 1\leq R$.

Set 
$C=\mathscr S/I_m^\times\mathscr S$.
Note, the natural map $A_v\to\mathscr S$ is flat. 
Since flatness preserves finite intersections of ideals, we have $
I_m^\times\mathscr S
=
\bigcap_{i=1}^{m}\ker(\text{res}_{v,e_i})\mathscr S$ and this gives the injection
\[
C\hookrightarrow\prod_{i=1}^{m}A_i^\an
\]
In particular, $C$ is reduced. The affinoid space associated to $C$ is the union
of these $m$ closed analytic branches. Its spectral
norm is therefore $
\|c\|_{\mathrm{sp}}
=
\max_{1\leq i\leq m}\|\text{res}_{v,e_i}(c)\|_{A_i}$.
Since $C$ is a reduced affinoid algebra, its spectral
norm is equivalent to its affinoid quotient norm.
Thus $C$ is complete for $\|\cdot\|_{ \mathrm{sp}}$.

The image of $A_v$ is dense in $\mathscr S$, so the image
of $A_v/\mathcal{T}_v$ is dense in $C$. Furthermore, the restriction
of $\|\cdot\|_{\mathrm{sp}}$ to $A_v/\mathcal{T}_v$ is exactly
$\|\cdot\|_{{S_v}}$. Hence the separated completion of $A_v$ with respect to $||\cdot||_{S_v}$ is isomorphic to $C$. 
Since $C$ is affinoid, it is Noetherian. Also, $
C\cong A_v/\mathcal{T}_v\otimes_{A_v}\mathscr S,$
so $A_v/\mathcal{T}_v\to C$ is flat.

It remains to prove that $C$ is Gorenstein. 
Let $X=\operatorname{Spec}(A_v/\mathcal{T}_v)$. The 
identification $C\cong A_v/\mathcal{T}_v\otimes_{A_v}\mathscr S$ realizes the completion $\widehat A_v^{\,\|\cdot\|_{S_v}}$  as
$\mathcal O(U)$, where
\[
U=\{x\in X^{\mathrm{an}}:
r\leq |x_j(x)|\leq R \text{ for }1\leq j\leq n\}
\]
is an affinoid domain in $X^{\mathrm{an}}$.
Since $X$ is Gorenstein, its analytification and every
affinoid domain therein are Gorenstein by
\cite[Lemma~2.4.6(1), Lemma--Definition~2.4.3(1)]
{ducros2018families}. Therefore, 
$\widehat{A}_v^{||\cdot||_{S_v}} \cong\mathcal O(U)$ is Gorenstein.
\end{proof}

\begin{rmk}
In Proposition \ref{gorenstein}, the branches need not have the same length.  More generally,
choose annuli with radii $
0<r_i\leq1\leq R_i, r_i,R_i\in|K^\times|$.
These define a Weierstrass domain $U$ inside the annulus
domain with radii $\min_i r_i$ and $\max_i R_i$ and the proof proceeds as usual. 
\end{rmk}

\begin{rmk}
    By the description of the kernel in Lemma \ref{kernel} and the proof of Proposition \ref{gorenstein}, the local model around a vertex of valence $d\geq 3$ is precisely the elliptic $m$-fold singularity of \cite[Definition 2.1]{smyth}. Following the idea that $\rho_*(\Oc_{C^\an})$ recovers the algebraic curve $C$, an interesting direction to pursue is to see if one can construct a proper singular analytic curve, hence algebrizable, from $\Oc_\sqsubset^\an$. This will be explored in future work. 
\end{rmk}

\begin{lem}
Let $A_v^\an$ denote $\Oc_\sqsubset^\an$ on the open star of $v$. Then the restriction morphism \begin{align*}
        A_v^\an \longrightarrow \oplus_i A_i^\an
    \end{align*}
    is injective. 
    \label{cechinjective}
\end{lem}

\begin{proof}
In the proof of Proposition \ref{gorenstein}, we proved injectivity for any closed star and thus by left exactness of inverse limits, we are done. 
\end{proof}

\begin{lem}
    There is a natural isomorphism $H^1(\Oc_\sqsubset^\times) \xrightarrow{\sim} H^1((\Oc_{\sqsubset}^\an)^\times)$. 
\end{lem}

\begin{proof}
    There is a short exact sequence 
    \begin{align*}
        1 \rightarrow \Oc_{\sqsubset}^\times \longrightarrow (\Oc_{\sqsubset}^\an)^\times \longrightarrow Q \longrightarrow 1.
    \end{align*}
    We first prove that $H^0(Q) = 1$. Let $f$ be a global section of $(\Oc_\sqsubset^\an)^\times /\Oc_\sqsubset^\times$. Then on edges, $f$ has a representative which is a 1-unit \ie of the form $\alpha(1+r)$ where $|r|<1$. It will be sufficent to prove our claim for a single simple cycle of $\underline{\sqsubset}$. Restricting to a simple cycle, would give a global section on this simple cycle. From the genus one case, we know there are no such global 1-units and hence $f=1$.

    To prove $H^1(Q)=$ is trivial, note that by Lemma \ref{cechinjective}, $A_v^\an $ embeds into $(\oplus_i A_{e_i}^\an)^\times = \oplus_i(A_{e_i}^\an)^\times$ and is surjective onto each summand. Thus the \v{C}ech differential for $Q$, with respect to the open cover given by open stars whose intersections are the interiors of edges, is surjective. This gives $H^1(Q) = 1$, completing the proof. 
    \end{proof}

    \subsubsection{Symmetric powers}


\begin{defn}
    We denote the degree map $\deg : \Pic(\underline{\sqsubset},\Oc_\sqsubset^\an) \longrightarrow \Z$ and the degree $d$ component of $\Pic(\sqsubset,\Oc_\sqsubset^\an)$ by $\Pic^d(\sqsubset)$. 
\end{defn}In the following definition, we endow $\TroPic^d(\sqsubset)$ with the structure of a ringed space. We first do the tropical case.

\begin{construct}[Tropical case]

Given a tropical curve $\underline{\sqsubset}$, consider the canonical torsor action of $\text{TroJac}(\sqsubset)$ given by $\tau: \text{TroJac}(\sqsubset) \times \text{TroPic}^d(\sqsubset) \longrightarrow \text{TroPic}^d(\sqsubset), (P,Q) \mapsto P+Q$. For every \(x\in \TroPic^d(\sqsubset)\), there is a canonical map $$\tau_x:\TroJac(\sqsubset)\xrightarrow{\sim}\TroPic^d(\sqsubset), a\longmapsto x+a.$$ Given $U\subseteq \TroPic^d(\sqsubset)$, define a continuous function \(f:U\to\mathbb R\) to be integral affine precisely when $f\circ\tau_x:\tau_x^{-1}(U)\longrightarrow\mathbb{R}$ is integral affine on \(\text{TroJac}(\sqsubset)\), for every \(x\in \text{TroPic}^d(\sqsubset)\). This defines a sheaf  $\Aff_{\Z, \TroPic^d(\sqsubset)}$ of integral affine functions on \mbox{$\TroPic^d(\sqsubset)$} which as in the case of the tropical Jacobian can be lifted to a $K$-affine structure on $\TroPic^d(\sqsubset)$. The same construction as for the tropical Jacobian then also endows $\TroJac^d(\sqsubset)$ with the structure of a ringed space, which we will denote $\Pic^d(\sqsubset)$. Furthermore, the morphism $\text{Sym}^d\Pic^1(\sqsubset) \longrightarrow \Pic^d(\sqsubset)$ is a morphism of ringed spaces as it is completely determined by the respective morphism on the integral affine structures. 

\label{construct:Picdtrop}
\end{construct}

\begin{construct}[Algebraic case]
    \label{construct:Picdalg}
By \cite[Theorem 4.2]{brandtulirsch}, there is a commutative diagram \begin{equation*}
\begin{tikzcd}
      C_d^\an \ar[d] \ar[r] & \Pic^d(C^\an) \ar[d]\\ 
    \underline{\sqsubset}_d \ar[r]& \TroPic^d(\sqsubset).
\end{tikzcd}
\end{equation*}
The degree $d$-Picard group $\Pic^d(C)$ has a canonical scheme structure, defined by representing the degree-\(d\) Picard functor. This defines the analytic structure on $\Pic^d(C^\an)$. Then $\TroPic^d(\sqsubset)$ is endowed with the pushforward of the sheaf of analytic functions $\Oc_{\Pic^d(C^\an)}$ via the retraction map. The symmetric product $\underline{\sqsubset}_d$ is equipped with $\rho_*(\Oc_{C_d^\an})$ and thus the diagram above is a commutative diagram of ringed spaces. We can also equip the the $d$-th symmetric product of $\underline{\sqsubset}_{\text{nb}}$ with the pushforward of $\Oc_{C_d}^\an$ under the retraction and contraction map and have a similar commutative diagram. 
\end{construct}

We next prove that the map  $\Oc_{\Pic^1(\sqsubset)} \rightarrow \alpha_*^\trop \Oc_\sqsubset  \rightarrow \Oc_\sqsubset^\an$ extends to  $\Oc_{\Pic^1(\sqsubset)}^\an := \rho_*(\Oc_{\Pic^1(C)^\an})$, where recall we abuse notation and write $\Pic^1(C)$ for the algebraic variety defined by the period matrix over $K$ as in the case of the Jacobian.  
Let $E_m = \spec R_m^\times$ and consider the natural map $f: E_m \hookrightarrow \G_m^{m-1} \hookrightarrow \G_m^g$ which induces a map on analytifications $f^\an: E_m^\an \hookrightarrow (\G_m^g)^\an$. Let $U$ be the affinoid subdomain of $E_m^\an$ such that $\Oc(U) = \Oc_{S_v}^\an$ and $V$ an  affinoid subdomain of $(\G_m^g)^\an$ such that $f^\an(U) \subseteq V$. In particular, this gives a bounded morphism  $\Oc(V) \rightarrow \Oc(U)$. We thus have the following result.  

\begin{prop}
    The morphism $\alpha^\trop$ factors as 
\begin{equation}
\begin{tikzcd}
(\underline{\sqsubset},\Oc_\sqsubset^\an) \ar[r, "\alpha^\trop"]\ar[rd, "\alpha^\trop"] & (\Pic^1(\sqsubset), \Oc_{\Pic^1(\sqsubset)}) \\
& (\Pic^1(\sqsubset), \Oc_{\Pic^1(\sqsubset)}^\an)\ar[u],
\end{tikzcd}
\end{equation}
where we continue to denote the map $(\underline{\sqsubset},\Oc_\sqsubset^\an) \longrightarrow (\Pic^1(\sqsubset), \Oc_{\Pic^1(\sqsubset)}^\an)$ by $\alpha^\trop$. 

\end{prop}

\subsubsection{Dualizing complex}
We end this section by studying the dualizing complex. Let
\[
     p\colon
    (\underline{\sqsubset},{\Oc}_{\sqsubset}^\an)
    \longrightarrow
    \operatorname{Spec}K
\]
denote the structure morphism. 
The exceptional inverse
image $ p^{\,!}$ exists in the sense of
\cite[Theorem B]{spaltenstein}. The theorem below and its proof apply both to the tropical and algebraic case. In the algebraic case, the proof passes through for $\sqsubset$ and $\sqsubset_{\text{nb}}$. 

\begin{thm}\label{thm:modified-dualizing-complex}
The dualizing complex
\[
    \omega_{{\sqsubset^\an}}^\bullet
    :=
 p^{!}K
\]
is concentrated in a single degree. More precisely, there is an
$\Oc_{\sqsubset}^\an$-module
$\omega_{{\sqsubset^\an}}$ such that
\begin{equation}\label{eq:modified-dualizing-sheaf}
 p^{!}K
    \simeq
    \omega_{{\sqsubset^\an}}[1].
\end{equation}
\end{thm}

\begin{proof}
Let $U$ be a sufficiently small open star centered at a vertex $v$,
and let $e_1,\ldots,e_d$ be its branches. We write $U = \cup_{i=1}^d [v,x_i)$. Compactly
supported cohomology of ${\Oc}^\an_{\sqsubset}|_U$ is computed by
the two-term cellular complex
\begin{equation}\label{eq:compact-support-cellular-complex}
 A_v^\an
    \longrightarrow
    \bigoplus_{i=1}^d (A_{e_i}^\an)_{x_i},
\end{equation}
placed in degrees $0$ and $1$, with differential given by the edge
restriction maps. Suppose the restriction of the function is zero at $x_i$;  by definition of the stalk, it is zero on some neighborhood and thus by the identity theorem \cite[Proposition 4.2]{ardakov-ben-bassat}, the function is identically zero on the branch. Then Lemma \ref{cechinjective} gives that the original function is identically zero on the star. 
This implies
\[
    H_c^0(U,{\Oc}_{\sqsubset}^\an)=0.
\]
Since the underlying space is one-dimensional, one also has
$H_c^i(U,{\Oc}_{\sqsubset}^\an)=0$ for $i>1$. Consequently,
$R\Gamma_c(U,{\Oc}_{\sqsubset}^\an)$ is concentrated in degree
one. The same assertion holds on sufficiently small open intervals.

Let $ p_U\colon
(U,{\Oc}_{\sqsubset}^\an|_U)\to\spec K$ be the
structure morphism. By the defining adjunction for
$ p_U^{\,!}$,
\begin{align*}
    R\Gamma\left(U, p^{!}K\right)
    &\simeq
    R\operatorname{Hom}_{{\Oc}^\an_U}
    \left({\Oc}^\an_U, p_U^{!}K\right)\\
    &\simeq
    R\operatorname{Hom}_K
    \left(R\Gamma_c(U,{\Oc}^\an_U),K\right).
\end{align*}
Since $K$ is a field and
$R\Gamma_c(U,{\Oc}_U^\an)$ is concentrated in degree $1$, the
right-hand side is concentrated in the corresponding dual degree.
Taking the filtered colimit over sufficiently small neighborhoods of
each point shows that the cohomology sheaves of
$ p_U^{\,!}K$ vanish in every degree except one. Defining
\[
    \omega_{{\sqsubset}^\an}
    :=
    \mathcal H^{-1}\left(p^{!}K\right)
\]
with the standard cohomological convention gives
\eqref{eq:modified-dualizing-sheaf}.
\end{proof}
Given that the dualizing complex is concentrated in a single degree, we make the following definition which applies in both the algebraic and tropical case. 

\begin{defn}
    The \textit{dualizing sheaf} is defined as $\omega_{\sqsubset^\an}[1] := p^!K$. 
\end{defn}

In \S\ref{serreduality}, we shall prove that the sheaf $\omega_{\sqsubset^\an}$ is a line bundle on $\sqsubset^\an$.

\subsection{$K$-tropical Divisors}
\label{ktropdivisor}
In this section, we want to construct a notion of a divisor on $(\underline{\sqsubset}, \Oc_\sqsubset)$ and on $(\underline{\sqsubset}, \Oc_\sqsubset^\an)$. They will not agree as divisors except on points whose topological support is at 2-valent vertices. Nevertheless, in \S\ref{divisorlinebundle} we will show that they do agree up to linear equivalence. We will first discuss the tropical case. 
\subsubsection{Realizable points - the tropical case}

We denote the $K$-points of $\Jac(\sqsubset)$ by $\Jac(\sqsubset)(K)$. By Lemma \ref{Kpointsjac}, the $K$-points of $\Jac(\sqsubset)$
    is isomorphic to $(N_\R\times T(K))/L$. Given a $K$-tropical curve $\sqsubset$, the $K$-tropical Jacobian $\Jac(\sqsubset)$ is a ringed space but also has a realization as an abelian variety over $K$; we  denote the abelian variety by $\Jac(C)$. Forgetting, the group structure of $\Jac(C)$ we can consider the retraction map $$\rho : \Pic^1(C)^\an \longrightarrow \Pic^1(\sqsubset)$$
    where $\Pic^1(C)^\an$ is the same underlying variety as $\Jac(C)^\an$. In Theorem \ref{morphismringedabelian}, this morphism is proven to be of ringed spaces and thus induces a morphism on $K$-valued points $ \Pic^1(C)^\an(K) \longrightarrow \Pic^1(\sqsubset)(K)$. 
    \begin{lem}
        The induced morphism $\rho_J: \Pic^1(C)^\an(K) \longrightarrow \Pic^1(\sqsubset)(K)$ is injective. 
    \end{lem}
    
    \begin{proof}
   It is enough to check this for the retraction map on Jacobians $ \Jac(C)^\an \longrightarrow \Jac(\sqsubset)$. The map on $K$-valued points is $T(K)/L \longrightarrow  (N_\R\times T(K))/L$ given by $z \mapsto (\rho_J(z), z)$ and thus if $[(\rho_J(x),x)] = [(\rho_J(x'), x')]$, then there exists $\ell \in L$ such that $x = \ell x'$ and thus $[x] =[x']$ in $T(K)/L$.  
    \end{proof}

Not every $K$-valued point of $\Pic^1(\sqsubset)(K)$ corresponds to a degree 1 line bundle. We thus make the following definition. 

\begin{defn}
    The \textit{realizable} $K$-points $\Pic^1(\sqsubset)(K)_\text{real}$ of $\Pic^1(\sqsubset)$ are given by the image $$\Pic^1(\sqsubset)(K)_\text{real}:= \rho(\Pic^1(C)^\an(K)),$$ 
    where $\Pic^1(C)^\an$ is the underlying variety of $\Jac(C)$. We can similarly define the \textit{realizable} $K$-points of $\Jac(\sqsubset)$, which we denote $\Jac(\sqsubset)(K)_\text{real}$. 
\end{defn}


In the algebraic case, a point $P \in C(K)$ correspond to the degree 1 line bundle $\Oc_C(P)$. In the $K$-tropical setting, not every point $P \in \sqsubset(K)$
corresponds to a line bundle. We thus single out the $K$-points which correspond to line bundles. 
\begin{defn}
    The \textit{realizable} $K$-points of $\sqsubset$ are $$ \sqsubset(K)_\text{real}:= \alpha^{-1}(\Pic^1(\sqsubset)(K)_\text{real}) \cap V(\sqsubset)$$
    where the intersection denotes that we only take the $K$-points whose topological support is at a vertex of $\sqsubset$. Given a vertex $v$ of $\sqsubset$, we denote $\sqsubset_v(K)_\text{real}$ to be the realizable $K$-points which have support at $v$.  
\end{defn}




\begin{defn}
    We shall denote the $K$-tropical Abel map restricted to the realizable points as $\alpha: \sqsubset(K)_\text{real}\longrightarrow\Pic^1(\sqsubset)(K)_\text{real}$. 
\end{defn}



\begin{rmk}
    The algebraic closure of $K$ is given by the field of Puiseux series $\overline{K} = \C\{\{t\} \} = \bigcup_{n \geq 1 }\C((t^\frac{1}{n}))$. At various points in the paper, we may need to take a finite extension of $K$ of the form $\C((t^{\frac{1}{n}}))$ to ensure the existence of rational points. Given a tropical curve, such an extension will only need to occur finitely many times, corresponding to the vertices. We will allow ourselves to make these finite extensions freely. 
    \label{extension}
\end{rmk}

The next two lemmas address the liftability to $\sqsubset(K)_\text{real}$ for points supported at vertices of $\underline{\sqsubset}$. The first lemma is an easy consequence of uniformization. 
\begin{lemma}
Let $v \in V(\sqsubset)$ and consider $\alpha(v) \in \text{TroJac}(\sqsubset)$. Then there exists $P_v \in \Jac(\sqsubset)(K)_\text{real}$ whose underlying point is $v$. 
\label{realizablejacobian}
\end{lemma}

\begin{lemma}
Every vertex of $\underline{\sqsubset}$ has a realizable
$K$-point lying above it.
\label{liftCpoint}
\end{lemma}

\begin{proof}
Note $\sqsubset(K)_v = \{ \epsilon: \Aff_{K,\sqsubset,v} \rightarrow K^\times \ | \ \epsilon|_{K^\times} = \text{id}\}$ and similarly for  $\Pic^1(\sqsubset)(K)$. 
A vertex $v$ of $\underline{\sqsubset}$ defines a morphism $\ev_v : \Aff_{\Z,\sqsubset,v} \longrightarrow \R$ given by evaluation at $v$. We then have a commutative diagram \begin{equation*}
\begin{tikzcd}
        \Aff_{\Z,\Jac(\sqsubset),\alpha(v)} \ar[r] \ar[rd] & \Aff_{\Z,\sqsubset,v} \ar[d] \\ 
        & \R.
\end{tikzcd}
\end{equation*}
Given that $\Jac(\sqsubset)$, is the lift of the tropical Jacobian by the section $a \mapsto t^a$, this diagram completely determines the same diagram for the $K$-affine structure. This then gives the commutative diagram:
\begin{equation*}
\begin{tikzcd}
        \Aff_{K,\TroPic^1(\sqsubset),\alpha(v)} \ar[r] \ar[d] & \Aff_{K,\sqsubset,v} \ar[d] \\ 
       K^\times \ar[r,"\text{id}"]  & K^\times.
\end{tikzcd}
\end{equation*}
By Lemma \ref{realizablejacobian}, the map $\Aff_{K,\TroPic^1(\sqsubset),\alpha(v)} \longrightarrow K^\times$ can be chosen to correspond to a realizable $K$-point of $\Pic^1(\sqsubset)$. Thus $\Aff_{K,\sqsubset,v} \rightarrow K^\times$ corresponds to a realizable $K$-point above $v$.



\end{proof}

\begin{rmk}
    The same argument also shows given any vertex of $\underline{\sqsubset}_d$, there exists a realizable $K$-point above it. 
    \label{extendtod}
\end{rmk}

\subsubsection{Principal divisors}
\label{subsec:realizable-principal-divisors}



We next describe the kernel of tropicalization on the realizable
points of the Jacobian. Let
\[
\trop_J\colon
\Jac(\sqsubset)(K)_{\mathrm{real}}
\longrightarrow
    \text{TropJac}(\sqsubset)
\]
denote the map which remembers the underlying tropical point, and
set
\[
\Jac(\sqsubset)(K)_{\mathrm{real},0}
:=
\ker(\trop_J).
\]
The next lemma describes the fiber of this tropicalization map. 

\begin{lemma}
\label{lem:unit-kernel-jacobian}
There is a canonical isomorphism
\[
\Jac(\sqsubset)(K)_{\mathrm{real},0}
\cong
\Hom\bigl(M,\Oc_K^\times\bigr).
\]

\end{lemma}

\begin{proof}
Under the uniformization
\[
\Jac(C)(K)=T(K)/L,
\qquad
T(K)=\Hom(M,K^\times),
\]
a class represented by \(z\in T(K)\) retracts to the origin if and
only if $\trop(z)\in\trop(L)$, where $\trop$ denotes the usual tropicalization. Since $\trop|_L\colon
L\longrightarrow\trop(L)$
is injective, there is a unique \(\ell\in L\) such that $\trop(z)=\trop(\ell).$
Replacing \(z\) by \(z\ell^{-1}\) does not change its class in
\(T(K)/L\), and gives $
\trop(z\ell^{-1})=0$. 
Thus every element of
    \(\Jac(\sqsubset)(K)_{\text{real},0}\) has a unique representative in $\Hom(M,\Oc_K^\times)$
which proves the assertion.
\end{proof}

The next lemma addresses how changing a realizable \(K\)-point
without changing its underlying topological point affects its Abel
image.

\begin{lemma}
\label{lem:unit-variation}
Let $v\in\underline{\sqsubset}$, $
x\in\sqsubset_v(K)_\text{real}$
and let $\psi\in
\Hom\bigl(\check{\Lambda}_{\sqsubset,v},\Oc_K^\times\bigr).$
Define \mbox{\(x\cdot\psi\in \ \sqsubset_v(K)\)} by
\begin{align}
\varepsilon_{x\cdot\psi}(a)
=
\varepsilon_x(a)\,\psi(\overline a),
\label{vary}
\end{align}
where $\overline a$ denotes the image of $a\in\Aff_{K,\sqsubset, v}$ in
\(\check\Lambda_{\sqsubset,v}\). Then \(x\cdot\psi\) is realizable. Furthermore, let
\[
d\alpha_v^\vee\colon
M\longrightarrow\check\Lambda_{\sqsubset,v}
\]
be the map on integral cotangent lattices induced by \(\alpha\). 
Under the identification of
Lemma~\ref{lem:unit-kernel-jacobian}, one has
\begin{equation}
\label{eq:unit-variation-formula}
\alpha(x\cdot\psi)-\alpha(x)
=
\psi\circ d\alpha_v^\vee
\in
\Hom(M,\Oc_K^\times).
\end{equation}

\end{lemma}

\begin{proof}
Since \(\psi\) takes values in \(\Oc_K^\times\), multiplication by
\(\psi\) does not change the valuation of any affine function. Let $\eta: \Aff_{K,\TroPic^1(\sqsubset),\alpha(v)} \longrightarrow K^\times$ correspond to the realizable $K$-point of $\Pic^1(\sqsubset)(K)_\text{real}$ given by the image of $x$ under $\alpha$. As in (\ref{vary}), we define $\eta_{ \psi\circ d\alpha_v^\vee}$ by  $\eta_{ \psi\circ d\alpha_v^\vee}(b) = \eta(b) \cdot \psi(d\alpha_v^\vee(\overline{b}))$ where $\overline{b}$ is the image of $b \in \Aff_{K,\TroPic^1(\sqsubset),\alpha(v)}$ in $M$. Then we have the commutative diagram \begin{equation*}
    \begin{tikzcd}
        \Aff_{K,\TroPic^1(\sqsubset),\alpha(v)} \ar[r] \ar[d, "\eta_{ \psi\circ d\alpha_v^\vee}"] & \Aff_{K, \sqsubset,v} \ar[d,"\epsilon_{x\cdot \psi}"] \\ K^\times  \ar[r,"\text{id}"] & K^\times. 
    \end{tikzcd}
\end{equation*} 
Therefore, by  Lemma \ref{lem:unit-kernel-jacobian}, \(x\cdot\psi\) is realizable. 

The points \(\alpha(x\cdot\psi)\) and \(\alpha(x)\) have the same
underlying tropical point, so their difference belongs to
\(\Jac(\sqsubset)(K)_{\mathrm{real},0}\). Let \(m\in M\), and choose a lift $b_m\in\Aff_{K,\Jac(\sqsubset),\alpha(v)}$ of $m$. Then
\[
\begin{aligned}
\bigl(\alpha(x\cdot\psi)-\alpha(x)\bigr)(m)
&=
\frac{
 \varepsilon_{x\cdot\psi}(\alpha^\#b_m)
}{
 \varepsilon_x(\alpha^\#b_m)
}\\
&=
\psi\bigl(\overline{\alpha^\#b_m}\bigr)\\
&=
\psi\bigl(d\alpha_v^\vee(m)\bigr).
\end{aligned}
\]
This proves \eqref{eq:unit-variation-formula}.
\end{proof}

\begin{lemma}
\label{lem:K-tropical-correction}
 For every $\eta\in\Jac^0(\sqsubset)(K)_{\mathrm{real},0}$
there exist realizable $K$-points $$
x_1,\ldots,x_g,x_1',\ldots,x_g'\in\sqsubset(K)_{\mathrm{real}}$$
such that $x_i$ and $x_i'$ have the same underlying topological point
for all $i$
and
\[
\eta
=
\sum_{i=1}^g
\bigl(\alpha(x_i')-\alpha(x_i)\bigr).
\]

\end{lemma}

\begin{proof}
Choose a maximal spanning tree $T$ of $\underline{\sqsubset}$, and let $
e_1,\ldots,e_g$
be the edges of $\underline{\sqsubset}\setminus T$. Orient each $e_i$, and let
$\gamma_1,\ldots,\gamma_g$
be the associated fundamental cycles, oriented so that the
coefficient of \(e_i\) in \(\gamma_i\) is \(1\). Then $\gamma_1,\ldots,\gamma_g$ is a basis of $M=H_1(\Gamma,\mathbb Z)$. 
For an oriented edge \(e\), let $
e^\vee\colon M\longrightarrow\mathbb Z$
denote the corresponding coedge. By construction,
$e_i^\vee(\gamma_j)=\delta_{ij},$
and hence $
e_1^\vee,\ldots,e_g^\vee$
is the dual basis of \(N\).

Let \(v_i\) be the initial vertex of \(e_i\), which we can always arrange to be a 2-valent vertex by further subdividing $\underline{\sqsubset}$. Recall the restriction map:
\[
\delta_i := \rho_{v_i,e_i}\colon\check{\Lambda}_{\sqsubset,v_i}\longrightarrow \check{\Lambda}_{\sqsubset,e_i} \cong \Z,
\qquad
\delta_i((a_h)_h)=a_{h_i}.
\]
where the isomorphism $\check{\Lambda}_{\sqsubset,e_i} \cong \Z$ follows from the choice of orientation on the edge $e_i$. 
The standard coedge description of the tropical Abel map gives
\begin{equation}
\label{eq:abel-differential-coedge}
\delta_i\circ d\alpha_{v_i}^\vee
=
e_i^\vee\colon M\longrightarrow\mathbb Z.
\end{equation}
vector represented by \(m\in M\) is the coefficient of \(e_i\) in
\(m\).

By Lemma~\ref{liftCpoint}, choose $
x_i\in\sqsubset_{v_i}(K)_{\mathrm{real}}$. 
Under Lemma~\ref{lem:unit-kernel-jacobian}, regard \(\eta\) as a
homomorphism $
\eta\colon M\longrightarrow \Oc_K^\times$.  
Set $
u_i:=\eta(\gamma_i)\in \Oc_K^\times$, 
and define \[
\psi_i\colon\check{\Lambda}_{\sqsubset,v_i}
\longrightarrow \Oc_K^\times,
\qquad
\psi_i(\lambda)=u_i^{\delta_i(\lambda)}.
\]
Put $x_i':=x_i\cdot\psi_i$. By Lemma \ref{lem:unit-variation}, the points $x_i$ and $x_i'$ are both realizable and lie above the
same vertex $v_i$. Moreover, by Lemma~\ref{lem:unit-variation} and \eqref{eq:abel-differential-coedge}, the difference
\(\alpha(x_i')-\alpha(x_i)\) corresponds to the character $
m\longmapsto u_i^{e_i^\vee(m)}.$
If $m=\sum_{i=1}^g a_i\gamma_i,$
then $e_i^\vee(m)=a_i$, and hence
\[
\begin{aligned}
\prod_{i=1}^g
\bigl(\alpha(x_i')-\alpha(x_i)\bigr)(m)
&=
\prod_{i=1}^g u_i^{a_i}\\
&=
\prod_{i=1}^g\eta(\gamma_i)^{a_i}\\
&=
\eta(m).
\end{aligned}
\]
Thus $
\eta
=
\sum_{i=1}^g
\bigl(\alpha(x_i')-\alpha(x_i)\bigr),$ as required.
\end{proof}

\begin{rmk}
\label{rmk:only2-valent}
	In the proof of Lemma \ref{lem:K-tropical-correction}, we can always ensure that the realizable $K$-points $x_i, x_i'\in \sqsubset(K)_\text{real}$ all have support given by 2-valent vertices. In the next section, we will show that given a divisor $D$ we can find a linearly equivalent divisor whose support only consists of 2-valent vertices. 
\end{rmk}

We can now use the realizable locus to define a notion of divisor on $\sqsubset$. This is similar to a notion of a Weil divisor; in the next section we will define a Cartier divisor for our ringed space. 

\begin{defn}
    A \textit{$K$-tropical realizable divisor on $\sqsubset$} is a formal sum of realizable $K$-points of $\sqsubset$. We denote the $K$-tropical realizable divisors by $\text{Div}(\sqsubset)$. One can define the degree in the usual way and we say two $K$-tropical divisors realizable $D,D'$ of the same degree $d$ are \textit{equivalent} if their image in $\Pic^d(\sqsubset)(K)_\text{real}$ is the same; we denote this by $D \sim D'$. Let the divisors of degree $d$ be denoted $\text{Div}^d(\sqsubset)$. We will only work with realizable divisors from now on and will refer to them as simply \textit{$K$-tropical divisors}. 
\end{defn}

\begin{defn}
	For a degree-zero $K$-tropical divisor $
D=\sum_x n_x\cdot x $, 
write $\alpha_*(D):=\sum_x n_x\cdot \alpha(x)\in\Jac(\sqsubset)(K)_{\mathrm{real}}.$ Define the \textit{principal $K$-tropical divisors} by

\[
\operatorname{Prin}(\sqsubset)
:=
\ker\left(
\alpha_*\colon
\operatorname{Div}^0(\sqsubset)
\longrightarrow
\Jac(\sqsubset)(K)_{\mathrm{real}}
\right).
\]
Let $\operatorname{TroPrin}(\sqsubset)$ be the group of principal tropical divisors on $\sqsubset$. 

\end{defn}
Since the following diagram commutes 
\begin{equation*}
    \begin{tikzcd}
        \text{Div}(\sqsubset) \ar[r, "\alpha"]\ar[d,"\trop"] & \Pic(\sqsubset)(K)_\text{real} \ar[d,"\trop_J"] \\ 
        \text{TroDiv}(\sqsubset) \ar[r, "\alpha^\trop"] & \text{TroPic}(\sqsubset)
    \end{tikzcd}
\end{equation*}
we have the equality 
\begin{equation}
\label{eq:divisor-Abel-compatibility}
\operatorname{trop}_J\bigl(\alpha_{,*}(D)\bigr)
=
\alpha^\trop_*
\bigl(\operatorname{trop}_*(D)\bigr)
\end{equation}
and that
\(\operatorname{trop}_*\) maps
\(\operatorname{Prin}(\sqsubset)\) to
\(\operatorname{TroPrin}(\sqsubset)\). Theorem 1.1 of \cite{jacobian} states that the pushforward morphism $\text{Prin}(C) \rightarrow \text{TropPrin}(\sqsubset)$ on principal divisors is surjective. They give an explicit construction of a lift and in the case of tropical curves, it consists adding lifts of 2-valent points. We extend their argument to the $K$-tropical setting. 
\begin{thm}
\label{principal-K-tropical-surjectivity}
The morphism
\[
\operatorname{trop}_*\colon
\operatorname{Prin}(\sqsubset)
\longrightarrow
\operatorname{TroPrin}(\sqsubset)
\]
is surjective. More precisely, let $
\overline D
=
\overline D_1-\overline D_2
\in\operatorname{TroPrin}(\sqsubset),
$
where \(\overline D_1\) and \(\overline D_2\) are effective tropical
divisors of degree \(n\). Then there exist effective \(K\)-tropical
divisors \(D_1,D_2\) of degree \(n+g\) such that
\[
D_1-D_2\in\operatorname{Prin}(\sqsubset)
,\text{ and }
\trop_*(D_1-D_2)=\overline D.
\]

\end{thm}

\begin{proof}
Write
\[
\overline D_j
=
\sum_{p\in V(\sqsubset)}n_{j,p} \cdot p ,
\qquad j=1,2.
\]
By Lemma \ref{liftCpoint}, for every point
\(p\) occurring in the support of \(\overline D_1\) or
\(\overline D_2\), choose a 
realizable, point $x_p\in\sqsubset_p(K)_{ \mathrm{real}}$ lying above \(p\). Define
\[
D_j^0
:=
\sum_p n_{j,p}\cdot x_p,
\qquad j=1,2.
\]
Then $\trop_*(D_j^0)=\overline D_j.$ Define 
\[
D^0:=D_1^0-D_2^0, \quad z:=\alpha_*(D^0)
\in\Jac(\sqsubset)(K)_{\mathrm{real}}.
\]
By \eqref{eq:divisor-Abel-compatibility}, $ \trop_J(z)
=
\alpha_*(\overline D).$
Since \(\overline D\) is a principal tropical divisor, we have $\alpha_*(\overline D)=0.$
Therefore $
z\in\Jac(\sqsubset)(K)_{\mathrm{real},0}.$

Applying Lemma~\ref{lem:K-tropical-correction} to \(-z\) gives 
realizable points $x_1,\ldots,x_g,x_1',\ldots,x_g'$ such that $x_i$ and $x_i'$ have the same topological support 
and $\sum_{i=1}^g\bigl(\alpha(x_i')-\alpha(x_i)\bigr)=-z.$
Define
\[
E_1:=\sum_{i=1}^g(x_i'),
\qquad
E_2:=\sum_{i=1}^g(x_i).
\]
Then $\trop_*(E_1-E_2)=0$ and $
\alpha_*(E_1-E_2)=-z$. Now set
\[
D_1:=D_1^0+E_1,
\qquad
D_2:=D_2^0+E_2.
\]
Then \(D_1\) and \(D_2\) are effective divisors of degree \(n+g\),
and
\[
\begin{aligned}
\alpha_*(D_1-D_2)
&=
\alpha_*(D^0)+\alpha_*(E_1-E_2) =
0.
\end{aligned}
\]
Hence $
D_1-D_2\in\operatorname{Prin}(\sqsubset)$ and 
\[
\begin{aligned}
\operatorname{trop}_*(D_1-D_2)
&=
\operatorname{trop}_*(D^0)
+
\operatorname{trop}_*(E_1-E_2)\\
&=
\overline D+0\\
&=
\overline D.
\end{aligned}
\]
\end{proof}

To end this section, we prove a $K$-tropical analogue of Jacobi inversion. The tropical analogue is \cite[Theorem 6.5]{mikhalkinzharkov} \ie for $d\geq g$, the tropical $d$-fold Abel map $\alpha_d^\trop: \underline{\sqsubset}_d \longrightarrow \text{TroPic}^d(\sqsubset)$ is a surjection. By Remark 9.3.15(c) in \cite{mikhalkin2009tropical}, $\alpha_d$ is a surjective tropical morphism. In particular, this implies that the map on structure sheaves \[\Aff_{\text{TroPic}^d(\sqsubset), \Z} \longrightarrow \Aff_{\sqsubset_d,\Z}\] is injective. This implies that the map on $K$-affine structures is injective and thus we have the following result. 

\begin{thm}[Jacobi inversion]
    When $d\geq g$, $\sqsubset_d(K)_\text{real} \longrightarrow \Pic^d(\sqsubset)(K)_\text{real}$ is a surjection. 
    \label{thm:jacobiinversion}
\end{thm}

\subsubsection{The algebraic case}

We end this section by defining divisors in the algebraic case. In the algebraic case, we will directly define a divisor on $\sqsubset^\an = (\sqsubset,\rho_*(\Oc_{C^\an}))$. 

\begin{defn}
    Let $\underline{\sqsubset}$ be the essential skeleton of an algebraic curve $C$. Then the \textit{realizable $K$-rational points} of $\sqsubset^\an$, which we denote $\sqsubset^\an(K)_\text{real}$, is equal to the image of $C^\an(K)$ under the retraction map $\rho: C^\an \longrightarrow \sqsubset^\an$. In particular, the topological support of the points in $\sqsubset^\an(K)_\text{real}$ consists of vertices. A \textit{realizable $K$-tropical divisor} is a formal sum of realizable $K$-points of $\sqsubset^\an$. As above, we can define the degree and linear equivalence. 
 \end{defn}

 We can also define the realizable $K$-rational points for the ringed space $
 \sqsubset_\text{nb}^\an:= (\underline{\sqsubset}_\text{nb}, c_*\Oc_\sqsubset^\an)$ as the image of $C^\an(K)$ under $c\circ \rho$. 

 \begin{lem}
     There is a bijection between the realizable $K$-rational points of $(\underline{\sqsubset},\Oc_\sqsubset^\an)$ and $(\underline{\sqsubset}_\text{nb},c_*(\Oc_\sqsubset^\an))$ which induces an equivalence between divisors up to linear equivalence. 
     \label{bridgelessbijection}
 \end{lem}

 We will primarily work with $\sqsubset^\an$ when studying its divisors. In \S\ref{globalsections}, we will remark how the results concerning global sections of line bundles extend to the contracted curve $\sqsubset_{\text{nb}}$.

\subsection{Divisor-line bundle correspondence}
\label{divisorlinebundle}

In this section, we establish an explicit divisor-line bundle correspondence.

\begin{defn}
   Given a realizable $K$-point $P \in \sqsubset(K)_\text{real}$, let $\Oc_\sqsubset(P) \in \Pic^1(\sqsubset)$ denote the corresponding line bundle defined via the morphism $\sqsubset(K)_\text{real} \longrightarrow \Pic^1(C)(K)$. Extending this to a $K$-tropical divisor $D$ using the group structure on $\Pic(\sqsubset^\an)$, let $\Oc_{\sqsubset^\an}(D)$ denote the corresponding line bundle.  
\end{defn}

The algebraic case is straightforward so we discuss it first.

\subsubsection{The algebraic case}

By definition of the divisor and linear equivalence, the divisor-line bundle correspondence follows from the correspondence on $C^\an$. 

\begin{cor}
    The assignment $D \mapsto \Oc_{\sqsubset^\an}(D)$ gives a degree preserving isomorphism \begin{align*}
        \frac{\text{Div}(\sqsubset)}{\sim} \xrightarrow{\ \ \  \sim \ \ \ } \text{Pic}(\sqsubset). 
    \end{align*}
\end{cor}

\subsubsection{The tropical case}




Recall, an effective Cartier divisor on a ringed space $X$ is the data of a pair $(\mathcal{I}, \iota)$ where $\mathcal{I}$ is an invertible $\Oc_X$-module, $\iota : \mathcal{I} \hookrightarrow \Oc_X$ is a monomorphism of $\Oc_X$-modules and for every $x\in X$ there exists an open neighbourhood $U$ of $x$, a trivialization $\mathcal{I}|_U\simeq\mathcal O_U$ under which $\iota_U$ is identified with multiplication by a section $f\in\Oc_X(U)$, such that $\Oc_U\xrightarrow{\cdot f}\mathcal O_U$ is a monomorphism of sheaves.

Furthermore, there is a 1-1 correspondence between effective Cartier divisors and line bundle section pairs for an arbitrary ringed space, where the section is represented by a non-zero divisor.  Indeed, from a divisor \((\mathcal I,\iota)\), set
$$ L=\mathcal I^\vee=\mathcal Hom_{\mathcal O_X}(\mathcal I,\mathcal O_X). $$
The inclusion $ \iota:\mathcal I\to\mathcal O_X $
is itself an element of
$\Hom(\mathcal I,\mathcal O_X) =H^0(X,\mathcal I^\vee), $
so it gives the canonical section $ s_D\in H^0(X,L). $
Thus $\mathcal O_X(D):=\mathcal I^\vee. $
Conversely, given an invertible sheaf \(L\) and a regular section
$ s:\mathcal O_X\to L, $ one can dualize it $ s^\vee:L^\vee\longrightarrow \mathcal O_X. $
Locally, if \(L\simeq\mathcal O_U\) and \(s\) is represented by \(f\), then
$ s^\vee:\mathcal O_U\xrightarrow{\cdot f}\mathcal O_U. $
Since \(f\) is a non-zero-divisor, this is injective. Hence
$ L^\vee\hookrightarrow\mathcal O_X $ is an effective Cartier divisor in exactly the sense above.
The two constructions are inverse. Starting with \((L,s)\), the associated divisor has ideal sheaf \(L^\vee\), so
$\mathcal O_X(D) =(L^\vee)^\vee \simeq L, $
since \(L\) is invertible. Under this identification, the canonical section of \(\mathcal O_X(D)\) is exactly \(s\).

Given an effective Cartier divisor $D$ on $\sqsubset^\an$, there is a natural embedding $\Oc_{\sqsubset^\an}(-D) \hookrightarrow \Oc_{\sqsubset^\an}$ and after tensoring by a line bundle $L$ gives $L(-D) \hookrightarrow L$. Taking global sections then gives the inclusion $H^0(L(-D)) \subseteq H^0(L)$. 


\begin{eg}
    \label{eg:2-valentcartier} 
Let $v$ be a 2-valent vertex of $\underline{\sqsubset}$ and $P \in \sqsubset_v^\an(K)_{\text{real}}$. Since $v$ is $2$-valent, $\Oc_{\sqsubset,v}^\an$ is isomorphic to annulus and is defined by the vanishing of $ f_P=x-a $ where $a \in K^\times$. 
Then \(x-a\) is a non-zero-divisor and a non-unit in \(A_v^\an\). Given an open interval containing $P$, the image of $f_P$ is an analytic unit and therefore the ideal generated by \(x-a\) has quotient supported precisely at \(v\). This defines an effective Cartier divisor. Let $\Lc_P$ denote the corresponding line bundle on $\sqsubset^\an$.

Now consider $P \in \sqsubset_v(K)_\text{real}$. Since we have  $$\sqsubset_v(K)_\text{real} = \alpha^{-1}(\Jac(\sqsubset)(K)_{\text{real},\alpha(v)})$$ where $\Jac(\sqsubset)(K)_{\text{real},\alpha(v)}$ is the fiber of $\alpha(v)$ under $\Jac(\sqsubset)(K)_\text{real}\longrightarrow \Jac(\sqsubset)$ and that $\Jac(\sqsubset)(K)_\text{real}$ is in bijection with $\Jac(C^\an)(K)$, there is a bijection of the 2-valent supported divisors on $\sqsubset$ and on $\sqsubset^\an$. 
\end{eg}

\begin{lem}
    There is a bijection of the 2-valent supported divisors on $\sqsubset$ and on $\sqsubset^\an$. 
\end{lem}

\begin{prop}
    Given a point $\tilde{P} \in \sqsubset(K)_\text{real}$ whose support is a $k$-valent vertex $P$ of $\underline{\sqsubset}$ for $k\geq 3$, there exists a linearly equivalent $K$-tropical divisor $D_P = P' +\sum_{i}(P_i-Q_i)$ whose support is only 2-valent vertices.
    \label{prop:linequiv}
\end{prop}

\begin{proof}
    By Lemma \ref{kvalent2valent}, there exists a tropical divisor $P' + \sum_i (P_i-Q_i) $ linearly equivalent to the support of $P$ where $P', P_i,Q_i$ are all 2-valent vertices. Let $\overline{D} = P_v -P' + \sum_i(P_i-Q_i)$ where $P_v$ is the support of $P$.  In the proof of Theorem \ref{principal-K-tropical-surjectivity}, one first chooses lifts of $P,P',P_i,Q_i$; for the vertex $P_v$, we take the lift $\tilde{P}$ and choose arbitrary lifts for $P', P_i$ and $Q_i$ which we denote $\tilde{P}', \tilde{P}_i$ and $\tilde{Q}_i$ respectively. Then as in Theorem \ref{principal-K-tropical-surjectivity}, define 
    \begin{equation*}
        D_1^0 = \tilde{P}+ \sum_i\tilde{Q}_i,\quad  D_2^0 = \tilde{P'} + \sum_i \tilde{Q}_i, \quad D^0 = D_1^0 - D_2^0, \quad z = \alpha(D^0)
    \end{equation*}
    By Remark \ref{rmk:only2-valent} and following the proof of Theorem \ref{principal-K-tropical-surjectivity}, there exist divisors $E_1$ and $E_2$ whose support only consists of 2-valent points and $\alpha_*(E_1-E_2) = -z$ and $\trop_*(E_1-E_2) = 0$. Defining $D_1 = D_1^0 + E_1$ and $ D_2 = D_2^0+E_2$ gives $D_1-D_2 \in \text{Prin}(\sqsubset)$. Thus $ \tilde{P} \sim \tilde{P}' + \sum_i(P_i-Q_i) + E_2 - E_1$ and the right hand side only consists of $K$-points which have support at 2-valent points by construction. 
\end{proof}

\begin{cor}
    Every line bundle on $\sqsubset$ is isomorphic to $\Oc_\sqsubset(P_1+\dots +P_a - Q_1-\cdots -Q_b)$ for some $a,b\geq 0$ and $P_i,Q_j \in \sqsubset(K)_\text{real}$ and their support consists of $2$-valent vertices. 
    \label{2-valentsum}
\end{cor}

\begin{lem}
   Let $v$ be a 2-valent vertex of $\sqsubset$ and $P \in \sqsubset_v(K)_{\text{real}}$. Then $\Lc_P\cong \Oc_{\sqsubset^\an}(P)$ where $\Lc_P$ is defined in Example \ref{eg:2-valentcartier}
   \label{cartierweil}
\end{lem}

\begin{proof}
    This is clear as the point $P$ on $\sqsubset(K)_\text{real}$ is defined by a single equation which is the same as in Example \ref{eg:2-valentcartier} and from the equivalence $H^1((\Oc_{\sqsubset}^\an)^\times) \cong \Pic(\sqsubset). $
    \end{proof}

Corollary \ref{2-valentsum} and Lemma \ref{cartierweil} give an assignment $D \mapsto \Oc_{\sqsubset^\an}(D)$ from $K$-tropical realizable divisors on $\sqsubset$ to $\Pic(\sqsubset^\an)$. By definition of linear equivalence of $K$-tropical realizable divisors and the bijection of the $2$-valent supported $K$-tropical realizable and Cartier divisors on $\sqsubset$ and $\sqsubset^\an$ respectively, we have the divisor-line bundle correspondence in the tropical case. 

\begin{cor}
    The assignment $D \mapsto \Oc_{\sqsubset^\an}(D)$ gives a degree preserving isomorphism \begin{align*}
        \frac{\text{Div}(\sqsubset)}{\sim} \xrightarrow{\ \ \  \sim \ \ \ } \text{Pic}(\sqsubset^\an). 
    \end{align*}
\end{cor}

\begin{rmk}
    It is important to note that the $K$-tropical realizable divisors are defined on $\sqsubset$ and not on $\sqsubset^\an$. Nevertheless, in the next section we will define the linear system of an effective $K$-tropical realizable divisor and prove that its rank matches with the rank of the corresponding line bundle on $\sqsubset^\an$. 
\end{rmk}

\subsection{Global sections}

\label{globalsections}
\subsubsection{The tropical case}
Now given an arbitrary divisor $D$ and 2-valent supported $K$-tropical realizable point $P \in \sqsubset(K)_\text{real}$, we have the isomorphism \begin{align*}
    \Lc_P\otimes \Oc_{\sqsubset^\an}(D)  \cong \Oc_{\sqsubset^\an}(P)\otimes \Oc_{\sqsubset^\an}(D)\cong \Oc_{\sqsubset^\an}(D+P). 
\end{align*}

Furthermore, given a point $P \in \sqsubset(K)_\text{real}$ supported at a 2-valent vertex $v$, by Example \ref{eg:2-valentcartier}, there is a short exact sequence 
\begin{equation}
    0 \rightarrow \Oc_{\sqsubset^\an}(-P) \rightarrow \Oc_{\sqsubset^\an} \rightarrow K|_P \rightarrow 0. 
    \label{inductionses}
\end{equation}
By twisting by $\Oc_{\sqsubset^\an}(D)$, we have the short exact sequence 
	\begin{align}
		0 \rightarrow \Oc_{\sqsubset^\an}(D-P) \rightarrow \Oc_{\sqsubset^\an}(D) \rightarrow \Oc_{\sqsubset^\an}(D)|_{P} \rightarrow 0
		\label{sesP}
	\end{align} 
	where $O_{\sqsubset^\an}(D)|_P$ is  supported only at the support of $P$. 
	
	\begin{eg}
	\label{eulerchar}
 Suppose  that $P \sim P' + \sum_{i=1}^n(P_i-Q_i)$ as in Lemma \ref{prop:linequiv} and that $P_j\neq Q_k$, $P' \neq P_j,Q_k $ for any $j,k$ where $P', P_i,Q_i \in \sqsubset(K)_\text{real}$ denote $K$-rational realizable points all which have support at $2$-valent points. 
    	Then we have a short exact sequence \begin{align*}
		0 \rightarrow \Oc_{\sqsubset^\an}\left(-P' -\sum_i P_i\right) \rightarrow \Oc_{\sqsubset^\an} \rightarrow K_{P'} \oplus_i K \rightarrow 0,
	\end{align*}
	and twisting by $\Oc_{\sqsubset^\an}\left(\sum_i Q_i\right)$ gives \begin{align*}
		0 \rightarrow \Oc_{\sqsubset^\an}\left(-P\right) \rightarrow \Oc_{\sqsubset^\an}\left(\sum_i Q_i\right)\rightarrow K_{P'} \oplus_i K \rightarrow 0. 
	\end{align*}
	Taking Euler characteristics gives $\chi(\Oc_{\sqsubset^\an}\left(\sum_iQ_i\right)) = \chi(\Oc_{\sqsubset^\an}(-P)) + n+1$. By induction and using the short exact sequence (\ref{inductionses}),  we have that $\chi(\Oc_{\sqsubset^\an}(\sum_i Q_i)) = n+ \chi(\Oc_{\sqsubset^\an})$ and thus \begin{equation}
\chi(\Oc_{\sqsubset^\an})=\chi(\Oc_{\sqsubset^\an}(-P)) + 1.
\label{eq:eulercharacteristic}
	\end{equation}
\end{eg}

\begin{prop}
    Let $D$ be a $K$-tropical realizable divisor. Then $\chi(\Oc_{\sqsubset^\an}(D))$ is finite. In particular, $H^i(\Oc_{\sqsubset^\an}(D))$ is finite for all $i$.  
    \label{finitecohom}
\end{prop}

\begin{proof}
    By Lemma \ref{prop:linequiv}, the divisor $D$ can be written as $\sum_i P_i - \sum_jQ_j$ where $P_i,Q_j$ are supported at 2-valent vertices. There is the short exact sequence 
    \begin{equation*}
        0 \rightarrow \Oc_{\sqsubset^\an}(-\sum_iQ_i) \rightarrow \Oc_{\sqsubset^\an}\rightarrow \oplus_j K_{Q_j} \rightarrow  0.
    \end{equation*}
Tensoring by $\Oc_{\sqsubset^\an}(\sum_iP_i)$ gives \begin{equation*}
    0 \rightarrow \Oc_{\sqsubset^\an}(D) \rightarrow \Oc_{\sqsubset^\an}(\sum_iP_i)\rightarrow \oplus_jK_{Q_j} \rightarrow 0. 
\end{equation*}  
and since $\Oc_{\sqsubset^\an}(\sum_iP_i)$ by Example \ref{eulerchar} and $K_{Q_j}$ have finite cohomology groups, we are done. 

\end{proof}

Recall that given a divisor $D$ on an algebraic curve $C$, the linear system \[|D| = \{ E \in \text{Div}(C) \: \ E\geq 0, E \sim D\}\] has rank $r_C(D) = \max \{n \in \Z_{\geq 0}: |D-E| \neq \varnothing\ \forall E \in \text{Div}^n_+(C(K))  \} = h^0(\Oc_C(D))-1$. The rank of a tropical line bundle is defined similarly. We  introduce the follow definition for the $K$-tropical curve $\sqsubset$. To simply notation, we will denote the tropical $d$-fold Abel map by $\alpha_d$. 

\begin{defn}
	Let $D$ be a degree $d\geq 0$ divisor on $\sqsubset$. The \textit{linear system} $|D|$ is defined as $\alpha_{d}^{-1}(\Oc_\sqsubset(D)) \subseteq \sqsubset^d(K)_\text{real}$. The \textit{rank} of $D$ is $-1$ if $|D|= \emptyset$ and otherwise,  \begin{align*}
		r_\sqsubset(D) := \max \{n \in \Z_{\geq 0}: |D-E| \neq \varnothing\ \forall E \in \text{Div}^n_+(\sqsubset(K)_\text{real})  \}.
	\end{align*}
\end{defn}

The next result shows that the $K$-tropical rank defined above behaves as the algebraic rank and further connects extends the divisor line bundle correspondence on $\sqsubset$.

\begin{lem}
Let $D$ be an effective realizable divisor on $\sqsubset$.  Then $r_\sqsubset(D) = h^0(L)-1$. 
\label{rankglobalsectionsequality}
\end{lem}

\begin{proof}
	We follow the proof method of \cite[Lemma 2.4]{specialbaker}. Let $P \in \sqsubset(K)_\text{real}$ be supported at a 2-valent vertex and taking long exact sequences of (\ref{sesP}) gives $H^0(D-P) = \ker(\text{ev}_P)$ where $\text{ev}_P : H^0(\Oc_{\sqsubset^\an}(D)) \rightarrow H^0(\Oc_{\sqsubset^\an}(D)|_P) \cong K$. 
If $h^0(\Oc_{\sqsubset^\an}(D)) \geq n+1$, it follows that for any points $P_1,\ldots,P_n \in \sqsubset(K)_\text{real}$ whose support is a 2-valent vertex, 
we have $ h^0(D-P_1-\cdots-P_n)\geq 1$, 
so that $H^0(D-P_1-\cdots-P_n)\neq (0)$ and
$|D-P_1-\cdots-P_k|\neq\varnothing$.
	
Conversely, we prove by induction on $n$ that if $h^0(L) = n$, then there exists $P_1,\dots, P_n \in \sqsubset(K)_\text{real}$ such that $h^0(L-\sum_iP_i) = (0)$ \ie $|D - \sum_i P_i | = \emptyset$. For $n = 0$, there is nothing to prove. Thus assume $h^0(L) = n \geq 1$ and choose a non-zero global section $s \in H^0(L)$. Dualizing the global section $s^\vee : L \rightarrow \Oc_{\sqsubset^\an}$, one can then either find or add a 2-valent vertex $v$ to $\underline{\sqsubset}$ with a realizable point $P$ in the fiber above $v$ such that $P^\sharp \circ s^\vee \neq 0$. The point $P$ defines an effective Cartier divisor and thus a degree 1 line bundle $\Oc_{\sqsubset^\an}(P)$ on $\sqsubset^\an$. In particular, $s \not\in H^0(\Oc_{\sqsubset^\an}(D-P_1))$ and there is a natural inclusion $H^0(\Oc_{\sqsubset^\an}(D-P_1)) \subseteq H^0(\Oc_{\sqsubset^\an}(D))$.  By induction, there exists $P_2,\dots, P_n$ such that $|D-P_1-\dots -P_n| = \emptyset$, thus completing the proof. 
\end{proof}

\begin{rmk}
   Note the lemma also proves that the cohomology of the line bundles on $\sqsubset_\text{nb}$ (contracting all bridges) is the same as the cohomology of the corresponding line bundle on $\sqsubset$. It is thus sufficient to work with $\sqsubset_\text{nb}$ for our purposes. 
    \label{workwithoutbriges}
\end{rmk}

\begin{lem}
If $D \in \text{Div}(\sqsubset)$, then \begin{align*}
		r_\sqsubset(D-P) \geq r_\sqsubset(D) -1
	\end{align*}
	for all realizable points $P \in \sqsubset(K)_\text{real}$
	and if $r_\sqsubset(D) \geq 0$, then there exists a realizable point $P$ such that $r_\sqsubset(D-P) = r(D) - 1$. 
	\label{rankdropbyone}
\end{lem}

\begin{proof}
	The proof is exactly as that of Lemma 2.6 in \cite{specialbaker}. 
\end{proof}

\subsubsection{The algebraic case}
The definitions and proofs of the results of the previous section easily pass through in the algebraic case. When needed, we will cite the $K$-tropical analogue. We will denote the rank $r_{\sqsubset}$ in the algebraic case. Baker's specialization lemma relates the rank of an algebraic line bundle to its tropicalization. We prove a similar result, first in the algebraic case, in the $K$-tropical setting.  The proof below is identical to the similar result in \cite{specialbaker}. We nevertheless include it as we will use the same argument in the tropical case in \S\ref{troprr}.

\begin{lem}[Specialization lemma]
	\label{special}
    	For all $D \in \text{Div}(C^\an)$, we have \[r_{\sqsubset}(\rho_*(D)) \geq r_{C^\an}(D).\]
\end{lem}
\begin{proof}
	  We prove by induction on $k$ that if $r_{C^\an}(D) = r_C(D) \geq k$, then $r_\sqsubset(\rho_*(D)) \geq k$. The base case $k=-1$ is clear so suppose that $k =0$ so that $r_C(D)\geq 0$. Thus there exists an effective divisor $E \in \text{Div}(C)$ such that $E \sim D$ so that $D-E \in \text{Prin}(C)$. Since $\rho_*$ is a a homomorphism and takes principal (resp. effective) divisors on $C$ to principal (resp. effective) divisors on $\sqsubset^\an$, we have $\rho_*(D) \sim \rho_*(E) \geq 0$ so $r_\sqsubset(\rho_*(D)) \geq 0$. 
	
	We may therefore assume that $k\geq 1$. Let $\overline{P} $ be a realizable point supported at a vertex $\rho(P) \in V(\sqsubset)$. By realizability, the image of $\overline{P}$ in $\Pic^1(\sqsubset)(K)_\text{real} = \Pic^1(C)(K)$ defines a degree 1 divisor  $P$ on $C$ such that $\rho_*(P) \sim \overline{P}$. Then $r_C(D-P) \geq k-1$, so by induction, we have $r_{\sqsubset^\an}(\rho_*(D) - \rho_*{P}) \geq k-1$. Since this is true for all realizable points $\overline{P} \in \sqsubset^\an(K)_\text{real},$ it follows from Lemma \ref{rankdropbyone} that $r_\sqsubset(\rho_*(D)) \geq k$ as required. 
	
\end{proof}

\begin{rmk}
   There is an alternative argument to prove the specialization lemma. Recall, we have the commutative diagram for every $d$:
\begin{equation}
	\begin{tikzcd}
		C_d^\an(K) \ar[r] \ar[d]  & \Pic^d(C^\an)(K) \ar[d]\\
		\sqsubset_{d}^\an(K)_\text{real}\ar[r] \ar[r] & \Pic^d(\sqsubset)(K)_\text{real}.
	\end{tikzcd}
    \label{isosquare}
\end{equation}
\label{easyspecialization}
Let $ D \in \text{Div}(C^\an)$, $E \in \text{Div}(\sqsubset)$ its image under the retraction map on $K$-valued points and suppose \(r=r_{C^{\mathrm{an}}}(D)\). For any effective divisor \(E'\) of degree \(r\) downstairs, choose an effective lift \(D'\) of degree \(r\). By the definition of rank, $
|D-D'|\neq\varnothing$. Specializing an effective representative gives $
|E-E'|\neq\varnothing.$ Since \(E'\) was arbitrary, $
r_{{\sqsubset}}(E)\ge r_{C^{\mathrm{an}}}(D).$
The case \(r=-1\) is automatic.

\end{rmk}

\begin{thm}
    Let $D \in \text{Div}(C^\an)$, $ E \in \text{Div}(\sqsubset)$ and suppose that $\rho_*(D) = E$. Writing $\Oc_{\sqsubset^\an}(E)$ and $\Oc_{C^\an}(D)$ for the corresponding line bundles on $\sqsubset$ and $C^\an$ respectively we have \begin{align*}
        h^0(\Oc_{\sqsubset^\an}(E)) = h^0(\Oc_{C^\an}(D)).    \end{align*}
        \label{thm:gagastyle}
\end{thm}

\begin{proof}
    This follows immediately from the commutative diagram (\ref{isosquare}) since the vertical arrows are bijections. 
\end{proof}

\begin{rmk}
    Note that Theorem \ref{thm:gagastyle} also applies to $\sqsubset^\an_\text{nb}$ since we still have the commutative diagram \ref{isosquare} in this case given the bijection between $K$-rational realizable points established in Lemma \ref{bridgelessbijection}. \end{rmk}

\subsubsection{Cheap $K$-tropical Riemann-Roch}
We first consider the tropical case.

\begin{prop}(Cheap $K$-tropical Riemann-Roch)
    Let $\Lc \in \Pic(\sqsubset^\an)$. Then \begin{equation}
        \chi(\sqsubset^\an, \Lc) = \deg \Lc + \chi(\underline{\sqsubset},\Oc_\sqsubset^\an)
    \end{equation}
    \label{cheapktroprr}
\end{prop}

\begin{proof}
Recall the short exact sequence (\ref{sesP}) 
\begin{align*}
0 \rightarrow \Lc(-P) \rightarrow \Lc \rightarrow K_P \rightarrow 0   
\end{align*}
where $P \in \sqsubset(K)_\text{real}$ is supported at a 2-valent vertex gives $\chi(\Lc) = \chi(\Lc(-P)) + 1$. By Proposition \ref{prop:linequiv}, we can write  
$$\Lc \cong \Oc(P_1+\cdots +P_a - Q_1-\dots -Q_b),$$ where $P_i,Q_j$ are supported at 2-valent vertices.
Iteratively applying $\chi(\Lc) = \chi(\Lc(-Q)) +1$ and $\chi(\Lc(P)) = \chi(\Lc) + 1$ gives 
\begin{align*}
    \chi(\Lc)= \chi(\Oc_\sqsubset^\an) + a-b = \chi(\Oc_\sqsubset^\an) + \deg \Lc. 
\end{align*}
\end{proof}

\begin{rmk}
    Given a degree one line bundle $\Lc \in \Pic^1(\sqsubset)(K)_\text{real}$ and choosing a base point $P_0 \in \sqsubset(K)_\text{real}$, consider the line bundle $\Lc' := \Lc \otimes \Oc_\sqsubset((g-1)P_0)$. This is a degree $g$ line bundle and thus by Jacobi inversion \ref{thm:jacobiinversion}, there exists $P_1, \dots, P_g \in \sqsubset(K)_\text{real}$ such that $\Lc'\cong  \Oc_\sqsubset(P_1+\dots +P_g)$. Thus $\Lc \cong \Oc_\sqsubset(P_1+\dots + P_g - (g-1)P_0)$. 
\end{rmk}

The cheap $K$-tropical Riemann-Roch theorem in the algebraic case also follows similarly.

\section{Serre Duality and $K$-tropical Riemann-Roch}
\label{serreduality}

In this section, we will prove $K$-tropical Serre duality and deduce the $K$-tropical Riemann-Roch from it. To do this, we first need to prove that the dualizing complex, which we proved is a sheaf concentrated in degree one in Theorem \ref{thm:modified-dualizing-complex}, is invertible. 
\subsection{Dualizing sheaf}
\label{dualizingsheaf}

One can compute the upper shriek functor via a closed embedding. In this section, we will assume we have contracted all the bridge edges in both the tropical and algebraic case. Recall, that the $K$-tropical Abel map factorises as:
    \begin{equation*}
         \begin{tikzcd}
        (\underline{\sqsubset},\Oc_\sqsubset^\an) \ar[r] \ar[rd] & (\Pic^1(\sqsubset), \Oc_{\Pic^1(\sqsubset)})
        \\ 
       & (\Pic^1(\sqsubset), \Oc_{\Pic^1(\sqsubset)}^\an)\ar[u].
    \end{tikzcd}
    \end{equation*}

We will then compute the upper shriek functor using $\alpha:(\underline{\sqsubset},\Oc_\sqsubset^\an)\longrightarrow (\Pic^1(\sqsubset), \Oc_{\Pic^1(\sqsubset)}^\an)$ and for ease of notation, will continue to denote it by $\alpha$. 

\begin{lem}
    The $K$-tropical Abel map is a closed immersion of ringed spaces.
\end{lem}

\begin{proof}
    Since we contracted all bridge edges in both the tropical and algebraic case, $\alpha$ is a closed immersion topologically. For the tropical case, by Construction \ref{modifedaffinestructure}, the map on rings is also surjective. The algebraic case follows since the analytic Abel map is a closed immersion. 
\end{proof}
Let $p_\sqsubset: \sqsubset^\an \longrightarrow \spec K$ denote the structure morphism. 
Hence by \cite[\href{https://stacks.math.columbia.edu/tag/0A74}{Tag 0A74}]{stacks-project}, we have \begin{align*}
   p_{\sqsubset}^!K = \alpha^{-1} \text{RHom}_{\Oc_{\Pic^1(\sqsubset)}^\an}\left(\alpha_*\Oc_\sqsubset^\an, p_{\Pic^1(\sqsubset)}^!K\right).
\end{align*}
     Then given the commutative diagram \begin{equation*}
         \begin{tikzcd}
             (\underline{\sqsubset},\Oc_{\sqsubset}^\an) \ar[r,"\alpha"] \ar[rd,"p_\sqsubset"] & (\Pic^1(\sqsubset), \Oc_{\Pic^1(\sqsubset)}^\an) \ar[d,"p_{\Pic^1(\sqsubset)}"] \\ 
             & \spec K
         \end{tikzcd}
     \end{equation*}
     and functoriality of the upper shriek functor, 
    \begin{align*}
        p^!_{\sqsubset} K&:= \alpha^!(p_{\Pic^1(\sqsubset)}^! K) 
    \end{align*}


One can also compute $p^!_{\Pic^1(\sqsubset)}K$ via a closed immersion as is done in the scheme-theoretic case. The computation for the $K$-tropical case is identical and thus we have the isomorphism \begin{align*}
    p^!_{\Pic^1(\sqsubset)}K  \cong \Oc_{\Pic^1(\sqsubset)}^\an[g],
\end{align*}
and thus there is an isomorphism \begin{align}
   p^!_\sqsubset K = \text{RHom}_{\Oc_{\Pic^1(\sqsubset)}^\an}\left(\alpha_*\Oc_\sqsubset^\an, \Oc_{\Pic^1(\sqsubset)}^\an[g]\right).
    \label{Rhom}
\end{align}






\begin{thm}
   The dualizing sheaf $p^!_{\sqsubset}K$ is equal to $\omega_{\sqsubset^\an}[1]$ for a line bundle $\omega_{\sqsubset^\an}$. 
\end{thm}

\begin{proof}
    By Theorem \ref{thm:modified-dualizing-complex}, $p^!_\sqsubset K = \omega_{\sqsubset^\an}[1]$ for some sheaf $\omega_{\sqsubset^\an}$. For the algebraic case, we know that $\alpha$ is a regular closed immersion and thus by \cite[\href{https://stacks.math.columbia.edu/tag/0BQZ}{Tag 0BQZ}]{stacks-project}, $\alpha^!$ sends a line bundle to a line bundle. By functoriality of the upper shriek functor, $p^!_\sqsubset K$ is a line bundle.  
    
    We next address the tropical case. It is sufficient to check this condition locally on closed sets. The result is clear for a closed interval strictly contained in an edge since it follows from the algebraic case. Thus we can assume we are working on a closed star $S_v$ of a vertex $v$ \ie $ p_v
    : (S_v,\Oc_{S_v}^\an)\longrightarrow \spec K $.   In Proposition \ref{gorenstein}, we proved that $\Oc_{S_v}^\an$ is Gorenstein. We compute the dualizing complex via a closed immersion as above.  Indeed, the natural map $A_v\to\mathscr S$ is flat, and
$\mathscr S$ is a regular Noetherian ring of dimension $m-1$ where $m$ is the valence of $v$. Let
\[
i_0:\operatorname{Spec}(A_v/\mathcal{T}_v)\hookrightarrow\operatorname{Spec}(A_v),
\qquad
i:\operatorname{Spec}(\Oc_{S_v}^\an)
\hookrightarrow\operatorname{Spec}(\mathscr S)
\]
be the closed immersions associated to the quotient maps where $i$ is the base change of the $i_0$ along
$A_v\to\mathscr S$. By flat base change, 
\cite[\href{https://stacks.math.columbia.edu/tag/0BZN}{Tag 0BZN}]{stacks-project} we have
\[
R\operatorname{Hom}_{\mathscr S}
(\Oc_{S_v}^\an,\mathscr S)[m-1]
\cong
\left(
 R\operatorname{Hom}_{A_v}(A_v/\mathcal{T}_v,A_v)[m-1]
\right)\otimes_{A_v/\mathcal{T}_v}^{\mathbf L}\Oc_{S_v}^\an.
\]
By Lemma \ref{kernel}, $A_v/\mathcal{T}_v$ is Gorenstein and thus 
\[
 R\operatorname{Hom}_{A_v}(A_v/\mathcal{T}_v,A_v)[m-1] \cong \omega_{A_v/\mathcal{T}_v}[1]
\]
where $\omega_{A_v/\mathcal{T}_v}$ is an invertible $A_v/\mathcal{T}_v$-module. Thus $$
p_v^! K = i^!(\mathscr S[m-1])
\simeq
(\omega_{A_v/\mathcal{T}_v}\otimes_{A_v/\mathcal{T}_v}\Oc_{S_v}^\an)[1],$$
and $\omega_{A_v/\mathcal{T}_v}\otimes_{A_v/\mathcal{T}_v}\Oc_{S_v}^\an$ is invertible. This completes the proof. 

\end{proof}


We next prove Serre duality in both the algebraic and tropical case.

\begin{thm}(Serre Duality)
     Let $\Lc$ be a line bundle on $\sqsubset^\an$. Then 
    \begin{align}
        H^i(\Lc)^\vee \cong H^{1-i}( \omega_{\sqsubset^\an}\otimes \Lc^{-1}).
    \end{align}
\end{thm}

\begin{proof}
   By Proposition \ref{finitecohom}, we have that the cohomology of a line bundle is finite. Then from Theorem \ref{grothendieck duality}, by taking global sections we have 
    \begin{align*}
        \Hom_{D(K)}(Rp_*(\Lc)[i], K) \cong \Hom_{D(\Oc_\sqsubset^\an)}(\Lc[i], p_\sqsubset^!K ).
    \end{align*}
    The left hand side is equal to $\Hom_K(H^i(\Lc), K)$ whilst the right hand side is $\text{Ext}^{-i}(\Lc,p_\sqsubset^!K)$. Since $p_\sqsubset^!K = \omega_{\sqsubset^\an}[1]$, this gives 
 \begin{align*}
 	\text{Ext}^{-i}(\Lc, p_\sqsubset^!K ) = \text{Ext}^{1-i}(\Lc, \omega_{\sqsubset^\an}) \cong H^{1-i}(\omega_{\sqsubset^\an}\otimes \Lc^{-1}).  
 \end{align*}
 \end{proof}

 \begin{defn}
    Under the divisor-line bundle correspondence, let $K_\sqsubset$ denote a divisor such that $\omega_{\sqsubset^\an} = \Oc(K_\sqsubset) $. 
\end{defn}

\subsection{$K$-tropical Riemann Roch}

In this section, we give the proof of the $K$-tropical Riemann-Roch theorem. We first give the proof in the tropical case. 
	\begin{lem}
		We have \begin{align*}
	\chi(\Oc_\sqsubset) = \chi(\Oc_\sqsubset^\an) = 1-g.
		\end{align*}
        \label{chi1-g}
	\end{lem}

\begin{proof}


We have already established that $H^1(\sqsubset,\Oc_\sqsubset^\times) = \Pic(\sqsubset)$. Let $\sqsubset_\epsilon = \sqsubset \times_K K[\epsilon]/(\epsilon^2)$ be the infinitesimal thickening of $\sqsubset$. On $\sqsubset$, there is surjective morphism of sheaves of rings $$\Oc_{\sqsubset_\epsilon} \rightarrow \Oc_\sqsubset, a+\epsilon b \mapsto a $$ which passing to units gives $\Oc_{\sqsubset_\epsilon}^\times \rightarrow \Oc_\sqsubset^\times $ whose kernel is precisely $1+ \epsilon \Oc_\sqsubset$. Thus there is a short exact sequence of abelian groups
\begin{align*}
	0 \rightarrow 1+ \epsilon \Oc_\sqsubset \rightarrow \Oc_{\sqsubset_\epsilon}^\times \rightarrow \Oc_\sqsubset^\times \rightarrow 0. 
\end{align*}
Moreover, $1+ \epsilon \Oc_\sqsubset \cong \Oc_\sqsubset$ via $1+ \epsilon b \mapsto b$ since $\epsilon^2 =0$. Taking long exact sequences gives \begin{align*}
	0 \rightarrow H^1(\Oc_\sqsubset) \rightarrow \Pic(\sqsubset_\epsilon) \rightarrow \Pic(\sqsubset) \rightarrow 0.
\end{align*}

Restricting to the degree 0 component, this short exact sequence coincides with its algebraic analogue concerning the deformations of $A:=\Pic(C)$, as an abelian variety, 
\begin{align*}
    0 \rightarrow T_0A \rightarrow  A_\epsilon \rightarrow  A \rightarrow 0.
\end{align*}
Since the tangent space of $A$, is $g$-dimensional,  we have that $h^1(\Oc_\sqsubset) = g$. The same argument proceeds identically for $\Oc_{\sqsubset}^\an$. This completes the proof. 
\end{proof}

\begin{prop}
    Let $\Lc$ be a degree 0 line bundle on $\sqsubset^\an$. Then \begin{align*}
        h^0(\Lc) = \begin{cases}
            1 \text{ if } \Lc \cong \Oc_\sqsubset^\an \\ 
            0 \text{ otherwise}
        \end{cases}
    \end{align*}
    \label{h0trivial}
\end{prop}

\begin{proof}
We can assume that $\Lc = \Oc(\sum_iP_i - \sum_jQ_j)$ where $P_i,Q_i$ are supported at 2-valent vertices. Suppose $ s\in H^0(\sqsubset^\an,\Lc)$.  We can associate the zero divisor $(s)_0$ consisting of realizable $K$-rational points whose support is a 2-valent vertex and where $s$ vanishes. By the correspondence between effective Cartier divisors and line-bundle section pairs, we have $\Lc\cong \Oc_\sqsubset^\an((s)_0)$ and thus $\deg (s)_0 =0$. This can only happen if $(s)_0 = 0$ \ie $s$ does not vanish. This implies that $s: \Oc\rightarrow \Lc$ is a trivialization and thus $\Lc \cong \Oc_\sqsubset^\an$. 
\end{proof}

For the algebraic case, Lemma \ref{chi1-g} follows from the analogous algebraic computation and Proposition \ref{h0trivial} follows from Theorem \ref{thm:gagastyle}. We can then prove the $K$-tropical Riemann-Roch theorem in both cases.

\begin{thm}($K$-tropical Riemann-Roch)
	Let $\Lc\in \Pic(\sqsubset^\an)$. Then \begin{align*}
		h^0(\Lc) - h^0(\omega_{\sqsubset^\an}\otimes \Lc^{-1}) = \deg(\Lc) - g +1.
	\end{align*}
	\label{ktroprr}
\end{thm}

\begin{proof}
We have \begin{align*}
    h^0(\Lc) - h^0(\omega_{\sqsubset^\an} \otimes \Lc^{-1}) &= h^0(\Lc) - h^1(\Lc)\\
    &= \chi(\Lc) \\ &= \deg \Lc +\chi(\Oc_\sqsubset^\an) \text{ by } (\ref{cheapktroprr}) \\ &= \deg \Lc - g+ 1 \text{ by } (\ref{chi1-g})
\end{align*}
\end{proof}

\begin{cor}
    The degree of the canonical divisor on $\sqsubset$ is $2g-2$. 
\end{cor}

\subsubsection{The algebraic Riemann-Roch theorem} We end this section by deducing the algebraic Riemann-Roch theorem for $C$ from the $K$-tropical analogue. We first verify compatibility between the canonical divisors. 

\begin{cor}
    Let $K_\sqsubset$ denote a divisor corresponding to $\omega_{\sqsubset^\an}$. Then $\rho_*(K_{C^\an}) \sim K_\sqsubset$. 
\end{cor}

\begin{proof}
    Applying the $K$-tropical Riemann-Roch theorem \ref{ktroprr} to $\rho_*(K_{C^\an})$ gives 
    \begin{align*}
        h^0(\rho_*K_{C^\an}) - h^0(K_\sqsubset - \rho_*K_{C^\an}) = g-1.
    \end{align*}
    By Theorem \ref{thm:gagastyle}, $r_C(K_{C^\an}) = r_C(K_C) = r_{\sqsubset}(\rho_*K_{C^\an}) = g-1$. The divisor $K_\sqsubset - \rho_*K_{C^\an}$ is degree 0 and thus by the algebraic analogue of Proposition \ref{h0trivial}, we have that $K_\sqsubset\sim \rho_*{K_{C^\an}}$ as required. 
\end{proof}


\begin{cor}[Algebraic Riemann-Roch]
Let $C$ be a genus $g$ maximally degenerate curve over $K$ and $\Lc \in \Pic(C)$. Then \[h^0(\Lc) - h^0(\omega_C\otimes \Lc^{-1}) = \deg \Lc - g + 1.\]
\end{cor}

\begin{proof}
Suppose $\Lc= \Oc_{C^\an}(D)$; then $\rho_*(K_C-D) \sim K_\sqsubset- \rho_*(D)$. Then by Theorem \ref{thm:gagastyle}, we have 
    \begin{align*}
        h^0(\Lc) = h^0(\Oc_{C^\an}(D)) = h^0(\Oc_{\sqsubset^\an}(\rho_*(D)), \end{align*}
        \begin{align*}h^0(\omega_{C^\an}\otimes \Lc^{-1}) = h^0(\Oc_{C^\an}(K_{C^\an}-D)) = h^0(\Oc_{\sqsubset^\an}(K_\sqsubset-\rho_*D)).
    \end{align*}
Thus applying the $K$-tropical Riemann-Roch theorem \ref{ktroprr}, gives the algebraic Riemann-Roch theorem. 
\end{proof}

\section{Tropical Riemann-Roch}
\label{troprr}
In the final section of this paper, we will show how to deduce the tropical Riemann-Roch theorem from the $K$-tropical Riemann-Roch theorem. Recall, we can assume that the tropical curve has no bridge edges.

\subsection{Lifting tropical to $K$-tropical}

We recall the construction of the $K$-affine structure for the tropical Riemann-Roch theorem.
 Let $s : \R \rightarrow K^\times, a \mapsto t^a$ denote the section of $v_K : \C((t^\R)) \rightarrow \R$. Since $\text{Ext}^1$ is covariant in the second entry, we have a homomorphism \[\text{Ext}^1(\check{\Lambda}_\sqsubset, \R) \longrightarrow\text{Ext}^1(\check{\Lambda}_\sqsubset,K^\times).\]

Denote the image of $\Aff_{\sqsubset,\Z}$ by $\Aff_{\sqsubset,K}$ and the $K$-tropical curve associated to $\Aff_{\sqsubset,K}$ by $(\underline{\sqsubset},\Oc_\sqsubset^\an)$. The section $s$ then defines a morphism $\Aff_{\Z,\sqsubset} \rightarrow \Aff_{K,\sqsubset}$ via the commutative diagram 
\begin{center}
	\begin{tikzcd}
		0 \ar[r] & K^\times \ar[r] \ar[d] & \Aff_{K,\sqsubset} \ar[r],\ar[d] & \check{\Lambda}_\sqsubset \ar[r] \ar[d,"\text{id}"] & 0 \\
		0\ar[r] & \R \arrow[u, out=20, in=340, "s"] \ar[r] & \Aff_{\Z,\sqsubset} 
		\ar[r] \arrow[u, out=20, in=340]&\check{\Lambda}_\sqsubset \ar[r] & 0
	\end{tikzcd}
\end{center}

Taking long exact sequences, gives the morphism $\Pic(\sqsubset) \rightarrow  \text{TroPic}(\sqsubset)$ and a section $\text{TroPic}(\sqsubset) \rightarrow \Pic(\sqsubset)$. The morphisms are compatible with the degree map $\deg: \Pic(\sqsubset) \rightarrow \Z$ and $\deg: \text{TroPic}(\sqsubset) \rightarrow \Z$ and thus for every degree $d$, there is a commutative diagram of sets
\begin{equation}
	\begin{tikzcd}
    		\sqsubset_d(K)_\text{real} \ar[r]  \ar[d, two heads, "\rho"] & \Pic^d(\sqsubset)(K)_\text{real}\ar[d, two heads]\\
		\sqsubset_d \ar[r] & \text{TroPic}^d(\sqsubset)
	\end{tikzcd}
\end{equation}

The proof of the specialization lemma in the algebraic case passes through exactly to achieve the $K$-tropical specialization lemma either via the proof given or the alternative proof given in Remark \ref{easyspecialization}.
\begin{lem}($K$-tropical specialization lemma)
\label{ktropspecial}
    For all $D \in \text{Div}(\sqsubset)$, we have \[r_{\sqsubset,\trop}(\rho_*(D)) \geq r_{\sqsubset}(D).\]
\end{lem}
A vertex $v$ of $\underline{\sqsubset}$ defines a morphism $\ev_v : \Aff_{\Z,\sqsubset,v} \longrightarrow \R$ given by evaluation at $v$ and similarly for $\alpha(v)$. We then have a commutative diagram \begin{equation*}
\begin{tikzcd}
        \Aff_{\Z,\Pic^1(\sqsubset),\alpha(v)} \ar[r] \ar[rd] & \Aff_{\Z,\sqsubset,v} \ar[d] \\ 
        & \R.
\end{tikzcd}
\end{equation*}
Given that $\Jac(\sqsubset)$ is the lift of the tropical Jacobian by the section $s$, this diagram completely determines the same diagram for the $K$-affine structure
\begin{equation*}
\begin{tikzcd}
        \Aff_{K,\Pic^1(\sqsubset),\alpha(v)} \ar[r] \ar[rd] & \Aff_{K,\sqsubset,v} \ar[d] \\ 
        & K^\times.
\end{tikzcd}
\end{equation*}
The same construction applies to the $d$-th symmetric product to give the commutative diagram 
\begin{equation}
\begin{tikzcd}
        \Aff_{K,\Pic^d(\sqsubset),\alpha(v)} \ar[r] \ar[rd] & \Aff_{K,\sqsubset_d,v} \ar[d] \\ 
        & K^\times.
\end{tikzcd}
\label{realizable}
\end{equation}
By the uniformization of $\Pic^d(\sqsubset)$ as an abelian variety, we have a commutative diagram 
\begin{equation}
    \begin{tikzcd}
        (K^\times)^g \ar[r] \ar[d] & (K^\times)^g/L \ar[d]\\
        \R^g \ar[r]  \arrow[u, out=110, in=250] & \R^g/v_K(L)  \arrow[u, out=110, in=250]
    \end{tikzcd}
\end{equation}
where $L$ is a full dimensional lattice in $(K^\times)^g$ 
and thus $\Aff_{K,\Pic^d(\sqsubset),\alpha(v)} \longrightarrow K^\times$ corresponds to a realizable $K$-point of $\Pic^d(\sqsubset)$.  
Then by the commutative diagram (\ref{realizable}), $\Aff_{K,\sqsubset_d,v} \rightarrow K^\times$ is realizable. 
This then gives the commutative diagram for every $d$:


\begin{equation}
	\begin{tikzcd}
    		\sqsubset_d(K)_\text{real} \ar[r]\arrow[d, two heads, out=290, in=65]    & \Pic^d(\sqsubset)\arrow[d, two heads, out=290, in=70] (K)_\text{real}\\
		\sqsubset_d \ar[r]\arrow[u, hook']& \text{TroPic}^d(\sqsubset)\arrow[u,hook']
	\end{tikzcd}
        \label{sectioncommute}
\end{equation}

\begin{prop}
	For every effective tropical divisor $D^\trop$ of degree $d$ on $\underline{\sqsubset}$, \begin{align*}
		r_{\sqsubset,\trop}(D^\trop) = r_\sqsubset(D)
	\end{align*}
	where $D$ is the image of $D^\trop$ under $\text{TroPic}^d(\sqsubset) \rightarrow \text{Pic}^d(\sqsubset)$.  
    \label{tropequal}
\end{prop}

\begin{proof}
    We establish the reverse inequality of the $K$-tropical specialization lemma. 
Suppose $r_{\sqsubset,\trop}(D) = r$. If $|D^\trop-D'^\trop| \neq \varnothing$ where $D'^\trop \in \text{TropDiv}(\sqsubset)$, then the corresponding linear system $|D-D'|$ of $\sqsubset$ is also non-empty by  the commutative diagram (\ref{sectioncommute}). In particular there exists an effective $K$-tropical realizable $D'$ of degree $r$ such that $|D-D'|\neq 0$. Thus we must have that $r = r_{\sqsubset,\trop}(D^\trop) \leq r_\sqsubset(D)$ by definition of the rank.  This implies that $r_{\sqsubset,\trop}(D^\trop) = r_\sqsubset(D)$. 
\end{proof}

\begin{cor}
    Under the map $\text{TroPic}^{2g-2}(\sqsubset) \longrightarrow \Pic^{2g-2}(\sqsubset)$, we have that $K_\sqsubset^\trop$ gets mapped to a divisor linearly equivalent to $K_\sqsubset$. 
\end{cor}

\begin{proof}
    Let $K_\sqsubset^\text{lift}$ denote the image of $K_\sqsubset^\trop$ under the map  $\text{TroPic}^{2g-2}(\sqsubset) \longrightarrow \Pic^{2g-2}(\sqsubset)$. Then by Proposition \ref{tropequal}, we have $r_{\sqsubset,\trop}(K_\sqsubset^\trop) = r_\sqsubset(K_\sqsubset^\text{lift}) = g-1$. Applying the $K$-tropical Riemann-Roch theorem \ref{ktroprr} to $K_\sqsubset^\text{lift}$ gives that $h^0(K_\sqsubset^\text{lift}) - h^0(K_\sqsubset- K_\sqsubset^\text{lift}) = g-1$. Thus as $K_\sqsubset-K_\sqsubset^\text{lift}$ is degree 0, by Lemma \ref{h0trivial}, we have that $K_\sqsubset-K_\sqsubset^\text{lift} \sim 0$. 
\end{proof}

Together with the $K$-tropical Riemann-Roch theorem (\ref{ktroprr}), this implies the tropical Riemann-Roch theorem. 

\begin{thm}
	Let $D^\trop$ be a tropical divisor on $\sqsubset$. Then 
	\begin{align*}
		r_{\sqsubset,\trop}(D^\trop) - r_{\sqsubset,\trop}(K_\sqsubset^\trop-D^\trop) = \deg D^\trop - g + 1. 
	\end{align*}
\end{thm}

\subsection{Riemann-Roch for finite graphs}

We end by deducing the Riemann-Roch theorem for finite graphs. Let $G$ be a finite graph without loop edges of genus $g$. Let $D$ be a divisor on $G$ and $r(D)$ denote the rank of $D$ on $G$. 
\begin{thm}\cite{luo,hladkykravnorine}
    Let $ \sqsubset$ be the tropical curve whose underlying graph is $G$ and whose edges all have length $1$. Then for every divisor $D$ on $G$, \begin{align*}
        r_{\sqsubset,\trop}(D) = r_G(D).
    \end{align*}
\end{thm}

Thus the tropical Riemann-Roch theorem applied to $\sqsubset$ as given in the theorem above implies the Riemann-Roch theorem for finite graphs. 

\begin{thm}
    Let $D$ be a divisor on $G$ and $K_G$ the canonical divisor of $G$. Then \begin{align*}
        r_G(D) - r_G(K_G-D) = \deg D -g +1
    \end{align*}
\end{thm}

\bibliographystyle{alpha}
\bibliography{zbib}

\end{document}